# Normalisation des opérateurs d'entrelacement et réductibilité des induites de cuspidales; le cas des groupes classiques $p$-adiques

Par C. Mœglin

Dans cet article, $F$ est un corps local non archimédien et on note $G(n)$ un groupe classique de rang $n$, défini sur $F$, c'est-à-dire le groupe symplectique $\mathrm{Sp}(2n)$ ou le groupe orthogonal d'une forme de dimension $2n$ ou $2n+1$ non nécessairement déployée (on remplacera le groupe orthogonal d'une forme de dimension impaire par sa composante connexe, le changement est anodin). On note $\pi_0$ une représentation cuspidale irréductible de $G(n,F)$. On fixe aussi $c \in \mathbb{N}$ et $\rho_0$ une représentation cuspidale unitaire irréductible de $\mathrm{GL}(c,F)$. On voit $\mathrm{GL}(c,F) \times G(n,F)$ comme le Levi d'un sous-groupe parabolique maximal de $G(n+c)$ où $G(n+c,F)$ est le groupe des automorphismes de la forme précédente à laquelle on a ajouté $c$ plans hyperboliques. Pour tout $s \in \mathbb{C}$, on définit l'induite $\rho_0|\det_{\mathrm{GL}(c,F)}|^s \times \pi_0$. Le problème est de déterminer les points $s$ où cette induite est réductible et d'en déduire une normalisation des opérateurs d'entrelacement standards:

$$\rho_0|\det_{\mathrm{GL}(c,F)}|^s \times \pi_0 \to \rho_0^*|\det_{\mathrm{GL}(c,F)}|^{-s} \times \pi_0,$$

normalisation notée $N(\rho_0 \times \pi_0, s)$ de telle sorte que cet opérateur soit holomorphe pour tout $s$ tel que $\mathrm{Re}\, s \geq 0$ et telle que le produit:

$$N(\rho_0^* \times \pi_0, -s)N(\rho_0 \times \pi_0, s) = 1.$$

C'est un programme initié par Langlands et repris par Shahidi. L'ennui est que même pour $\mathrm{SL}(2,F)$ il semble difficile d'obtenir la formule du produit sans faire des choix pas complètement naturels (remarque de Waldspurger). Dans l'article, on ne cherche pas à obtenir vraiment 1 dans la formule du produit mais à remplacer 1 par une fonction holomorphe inversible (c'est-à-dire sans zéros) que l'on ne sait pas calculer. Les travaux de Shahidi résolvent cette question si l'on suppose que $\pi_0$ a un modèle de Whittaker [S$_2$] et partiellement dans d'autres cas [S$_3$].

Cette article démontre effectivement les résultats ci-dessus mais avec l'hypothèse: que la situation locale provienne d'une situation globale satisfaisant une propriété de fonctorialité (forme faible de l'existence d'un relèvement à un groupe linéaire). Précisément, nous faisons l'hypothèse suivante:

On fixe $k$ un corps de nombres et $v_0$ une place de $k$ telle que $k_{v_0} \simeq F$; on note $\mathbb{A}$ les adèles de $k$. Pour $a, b \in \mathbb{N}$ et pour $\sigma$ une représentation cuspidale unitaire irréductible de $\mathrm{GL}(b, \mathbb{A})$ on note $\pi(\sigma, a)$ la représentation automorphe de carré intégrable de $\mathrm{GL}(ab, \mathbb{A})$, dite module de Speh, unique quotient irréductible de la représentation induite: $\times_{k\in[0,a-1]}\sigma||^{(a-1)/2-k}$. On prend la même définition si $\sigma$ est une série discrète de $\mathrm{GL}(b, F)$ (cela a un sens grâce à [Z]). On pose $\varepsilon = 0$ sauf si $G = \mathrm{Sp}(2n)$ où $\varepsilon = 1$.

*Soit $\pi_0$, comme ci-dessus, une représentation cuspidale irréductible de $G(n, F)$ alors il existe une représentation $\tilde{\pi}$ de $\mathrm{GL}(2n+\varepsilon, \mathbb{A})$ irréductible avec les propriétés suivantes*:

(A)· *$\tilde{\pi}$ est $G(n)$-discrète notion que l'on va définir par les 2 propriétés suivantes*:

·  *il existe un ensemble fini, $E$, de couples $(\sigma, b)$ où $\sigma$ est une représentation cuspidale irréductible et unitaire d'un groupe linéaire $\mathrm{GL}(n_\sigma, \mathbb{A})$ (ce qui définit $n_\sigma$) et où $b \in \mathbb{N}$ tel que $\tilde{\pi}$ soit l'induite des modules de Speh $\pi(\sigma, b)$, i.e.*:

(1) $$\tilde{\pi} \simeq \times_{(\sigma,b)\in E}\pi(\sigma, b);$$

·  *notons $\Psi(\tilde{\pi})$ l'ensemble des représentations $\pi$ automorphes irréductibles de carré intégrable de $G^0(n, \mathbb{A})$ ($G^O$ est le groupe $G$ sauf dans les cas des groupes orthogonaux paires où il représentation le groupe spécial orthogonal) telles que pour toute place $v$ de $k$, hors d'un ensemble fini de places, $S$, contenant les places archimédiennes et arbitrairement grand, $\pi_v$ corresponde à $\tilde{\pi}_v$ par l'inclusion du L-groupes de $G(n, k_v)$ dans $\mathrm{GL}(2n+\varepsilon, \mathbb{C})$. Alors cet ensemble est non vide*;

(B)· *$\Psi(\tilde{\pi})$ (défini ci-dessus) contient une représentation $\pi$ telle que $\pi_{v_0} \simeq \pi_0$ (si $G = \mathrm{O}(2n, F)$ on demande uniquement que $\pi_{v_0}$ est l'une des composantes de la restriction de $\pi_0$ à $\mathrm{SO}(2n, F)$).*

Le premier point que l'on démontre (avec cette généralité) (cf. §4 ci-dessous) est qu'il est possible d'associer à $\pi_0$ une représentation de $\mathrm{GL}(2n+\varepsilon, F)$, $\tilde{\pi}_0^L$ indépendante du choix de $\tilde{\pi}$ vérifiant (A) et (B); cette représentation $\tilde{\pi}_0^L$ détermine conjecturalement (cf. §6 ci-dessous) le paquet de Langlands auquel $\pi_0$ appartient et c'est grâce à elle que l'on peut avec uniquement les hypothèses (A) et (B) normaliser les opérateurs d'entrelacement et calculer les points de réductibilité. Cette représentation est autoduale. On va être plus précis dans la suite de l'introduction mais c'est le moment d'insister sur le fait que l'on considère réellement les groupes orthogonaux pairs et pas seulement leur composante connexe; c'est le seul moyen d'unifier les résultats.

Expliquons d'abord en utilisant le yoga d'Arthur, en quoi $\tilde{\pi}_0^L$ est lié au paquet de Langlands contenant $\pi_0$: si les conjectures d'Arthur locales sont loin d'être connues pour les groupes du type $G(n)$ par contre on a une bonne idée de ce que sont les paquets d'Arthur pour $\text{GL}(2n+\varepsilon, F)$. Ces paquets sont réduits à un élément que l'on peut décrire explicitement, toute représentation n'est pas dans un paquet mais les composantes locales des représentations du type (1) y sont. Donc $\tilde{\pi}$ étant fixée comme ci-dessus, on peut décrire le paquet auquel doit appartenir sa composante locale $\tilde{\pi}_{v_0}$ (cf. §6); ce paquet est donné par un homomorphisme:

$$\psi_{\tilde{\pi}}: W_F \times \text{SL}(2,\mathbb{C}) \times \text{SL}(2,\mathbb{C}) \to \text{GL}(2n+\varepsilon, \mathbb{C})$$

On montre que la restriction de cet homorphisme au produit de $W_F$ avec la diagonale de $\text{SL}(2,\mathbb{C}) \times \text{SL}(2,\mathbb{C})$ est indépendante du choix de $\tilde{\pi}$. Cela donne un homomorphisme:

(∗) $\qquad W_F \times \text{SL}(2,\mathbb{C}) \to \text{GL}(2n+\varepsilon, \mathbb{C}),$

que l'on voit comme un paramètre de Langlands. Et $\tilde{\pi}_0^L$ est la représentation de Langlands attachée à ce paquet. C'est grâce à cette représentation ou, ce qui revient au même, à sa duale, $\tilde{\pi}_0^A$ pour l'involution de Zelevinski que l'on décrit les points de réductibilité pour les induites $\rho_0||^s \times \pi_0$. Remarquons que l'on ne montre pas dans cet article que l'image du morphisme en (∗) se factorise par $^L G(n)$ vu comme sous-groupe de $\text{GL}(2n+\varepsilon, \mathbb{C})$ (cf. §6); la seule chose qui manque est une propriété de la correspondance de Langlands locale pour GL(n) (maintenant démontrée par Harris-Taylor, puis Henniart) à savoir que le morphisme de $W_F \to \text{GL}(c, \mathbb{C})$ associé à une représentation cuspidale $\sigma$ d'un groupe $\text{GL}(c, F)$, se factorise par $\text{Sp}(c, \mathbb{C})$ (resp. $\text{O}(c, \mathbb{C})$ si le facteur $L_{v_0}(\rho_0, \wedge^2, s)$ (resp. $L_{v_0}(\rho_0, \text{Sym}^2, s)$) a un pôle en $s=0$ - au sujet de ces fonctions $L_{v_0}$ qui interviennent naturellement mais de façon clé dans notre article, on renvoie à Shahidi et spécialement à [$S_4$]; de plus notre article est très proche de [$S_4$]. Mais si l'on admet cela, l'homomorphisme (∗) se factorise par $^L G(n)$ et est naturellement l'homomorphisme qui doit être attaché à la représentation cuspidale $\pi_0$ par les conjectures de Langlands (cf. §6).

Revenons aux seules hypothèses (A) et (B). On déduit de la construction de $\tilde{\pi}_0^L$ et $\tilde{\pi}_0^A$, deux normalisations des opérateurs d'entrelacement. Le spécialiste sait déjà qu'il n'y a pas qu'une normalisation même en lui imposant d'être situation locale d'une situation globale; en fait il doit y en avoir au moins autant que de paquets d'Arthur contenant $\pi_0$ et la différence entre toutes ces normalisations possibles vient des zéros situés dans $\mathbb{R}_{>0}$; la normalisation "minimale" n'a pas de zéros dans cet espace et est celle de Langlands. La normalisation obtenue avec $\tilde{\pi}_0^A$ a des zéros (pour au moins une représentation cuspidale $\rho_0$) exactement quand $\tilde{\pi}_0^A \not\simeq \tilde{\pi}_0^L$ et ceci est encore équivalent à ce que $\tilde{\pi}_0^L$ soit une induite de représentations cuspidales unitaires. Bien sûr on a l'holomorphie pour tout $s \in \mathbb{R}_{\geq 0}$ et l'unitarité sur l'axe imaginaire.

On démontre l'indépendance de $\tilde{\pi}_0^A$ et $\tilde{\pi}_0^L$ du choix de $\tilde{\pi}$ en étudiant les points de réductibilité pour les induites comme ci-dessus; cette étude se fait en globalisant la situation (ce qui nécessite les hypothèses (A) et (B)) et en utilisant l'équation fonctionnelle pour les opérateurs d'entrelacement attachés à un parabolique maximal et à une représentation cuspidale de ce parabolique.

Venons-en à une description plus explicite. Les travaux de Shahidi suggèrent une généralisation précise pour la description de ces points de réductibilité. Les premières propriétés auxquelles on s'attend, sont les suivantes: il y a au plus 1 point de réductibilité (sur l'axe réel) et il y a exactement 1 point de réductibilité (réel) quand $\rho_0 \simeq \rho_0^*$. Ce point de réductibilité est lié aux propriétés d'endoscopie satisfaites par $\pi_0$. Plutôt que de parler d'endoscopie, il est plus simple d'avoir en tête les paramétrisations de Langlands et d'Arthur qui traduisent ces propriétés en termes géométriques. Avec les seules hypothèses (A) et (B), on n'arrive pas à ce résultat; il peut y avoir plus d'un point de réductibilité; cela est décrit avec précision en Section 4 ci-dessous mais n'a pas vraiment d'intérêt en soi; le seul intérêt est que c'est ce calcul qui permet de définir $\tilde{\pi}_0^L$ indépendamment du choix de $\tilde{\pi}$.

Après que ce papier eut été écrit, Shahidi et Muic (via Tadic) m'ont indiqué que Silberger a démontré dans [Si] qu'il y a au plus un point de réductibilité sur l'axe réel. En tenant compte de cela et étant donné notre Section 4, cela prouve que l'hypothèse (C″) expliquée ci-dessous est conséquence de (A) et (B) et que la description que nous faisons "du" point de réductibilité résulte uniquement de (A) et (B).

Dans les paragraphes 2 et 3, on montre sous quelles hypothèses on peut obtenir ce résultat pour les points de réductibilité. Bien sûr ces hypothèses sont très certainement vérifiées. En 2, on impose l'hypothèse suivante:

(C). Il existe $\tilde{\pi}$ satisfaisant (A) et (B) et tel qu'il existe un ensemble fini $E_0$ de couples $(\sigma_0, a)$, où $\sigma_0$ est une représentation cuspidale irréducible et unitaire d'un groupe linéaire $\mathrm{GL}(n_{\sigma_0}, F)$ et où $a \in \mathbb{N}$ tel que:

$$\tilde{\pi}_0 := \tilde{\pi}_{v_0} \simeq \times_{(\sigma_0, a) \in E_0} \pi(\sigma_0, a).$$

Il est clair que (A) et (B) sont les premières propriétés que l'on peut obtenir avec la formule des traces d'Arthur: partant de $\pi_0$, il est naturel de penser qu'il existe $\pi$ une représentation cuspidale de $G(n, \mathbb{A})$ telle que $\pi_{v_0} \simeq \pi_0$. Par la comparaison de la formule des traces stabilisée d'Arthur pour $G(n)$ avec la formule des traces stable pour un produit semi-direct de $\mathrm{GL}(2n+\varepsilon)$ avec un automorphisme extérieur convenable dépendant de $G(n)$ (cf. [A]), on espère construire $\tilde{\pi}$ qui vérifiera (A) et (B).

La propriété (C) peut sembler plus étrange. Il faut la comprendre de la façon suivante: On part de $\pi_0$ cuspidale et on espère savoir lui associer un homomorphisme de Langlands:

$$\phi_{\pi_0} : W_F \times \mathrm{SL}(2,\mathbb{C}) \to {}^L G(n).$$

On prolonge cet homomorphisme en un homomorphisme de $W_F \times \mathrm{SL}(2,\mathbb{C}) \times \mathrm{SL}(2,\mathbb{C})$ dans ${}^L G(n)$ trivial sur la deuxième copie de $\mathrm{SL}(2,\mathbb{C})$ puis on définit $\psi_{\pi_0}$ en échangeant les 2 copies de $\mathrm{SL}(2,\mathbb{C})$, pour obtenir:

$$\psi_{\pi_0} : W_F \times \mathrm{SL}(2,\mathbb{C}) \times \mathrm{SL}(2,C) \to {}^L G(n)$$

homomorphisme trivial sur la première copie de $\mathrm{SL}(2,\mathbb{C})$ et on préfère remplacer ce facteur par $\mathrm{SU}(2,\mathbb{R})$ (ce changement est totalement anodin). S'il existe un bon groupe Tannakien de $k$, noté $L_k$, ce groupe doit contenir $W_F \times \mathrm{SU}(2,\mathbb{R})$ et l'homomorphisme $\psi_{\pi_0}$ doit se prolonger en un homomorphisme à la Arthur:

$$\Psi : L_k \times \mathrm{SL}(2,\mathbb{C}) \to {}^L G(n).$$

Ici le $L$-groupe est celui associé à $G(n,k)$ et non plus $G(n,F)$. La définition même de $L_k$ doit permettre la construction de $\tilde\pi$ associé à $\Psi$; il faut maintenant utiliser la formule des traces pour montrer que $\tilde\pi$ vérifie (A) et (B) mais avec cette façon de faire (C) est automatique, il résulte simplement du fait que $\psi_{\pi_0}$ est la restriction de $\Psi$ à $W_F \times \mathrm{SU}(2,\mathbb{R}) \times \mathrm{SL}(2,\mathbb{C})$. Dans cette présentation, on peut quand même éviter de faire appel à $L_k$ qui est un peu trop conjectural et elle est en principe atteignable avec la formule des traces d'Arthur (cf. [A]).

La condition (C) vient de la philosophie générale relative aux paquets d'Arthur mais ce n'est sans doute pas celle-là que l'on obtiendra en premier mais une condition plus technique qui est la suivante: on renvoie le lecteur à [A] pour la définition de l'automorphisme extérieur $\theta$ de $\mathrm{GL}(2n+\varepsilon, k)$ tel que $G(n,k)$ soit un groupe endoscopique elliptique pour le produit semi-direct $\mathrm{GL}(2n+\varepsilon, k) \times <\theta>$. La condition que l'on devrait pouvoir démontrer avec la comparaison des traces est que si $\tilde\pi$ satisfait (A) et (B) alors $\tilde\pi_{v_0}$ n'est pas une induite irréductible à partir d'un parabolique $\theta$-stable de $\mathrm{GL}(2n+\varepsilon, F)$, propriété que l'on peut résumer en disant que $\tilde\pi_{v_0}$ est $G(n)$ discrète; ceci venant évidemment du fait que $\pi_0$ est cuspidal, cela ne peut être vrai en toute place. La condition que nous avons noté (C″) en Section 3 est encore plus faible que celle-là mais suffit pour obtenir les propriétés cherchées pour les points de réductibilité des induites que l'on va maintenant expliquer.

Remarquons toutefois pour les sceptiques, que l'on peut (cf. §3 ci-dessous) remplacer la condition (C) par la condition (C′) à laquelle tout le monde croit,

(C′) Les représentations cuspidales intervant dans le support cuspidal de $\tilde\pi$ vérifient la conjecture de Ramanujan-Petersson à la place $v_0$; c'est-à-dire avec les notations de (1) que les $\sigma$ sont telles que $\sigma_{v_0}$ est tempérée.

Supposons donc que l'une des hypothèses (C), (C′) ou (C″) soient vérifiées et décrivons les points de réductibilité. On a besoin de quelques notations. On a défini ci-dessus la représentation $\tilde{\pi}_0^A$; dans le cas où (C) est vérifiée c'est précisément la représentation $\tilde{\pi}_{v_0}$ et dans tous les cas elle est de la forme:

$$\tilde{\pi}_0^A \simeq \times_{(\sigma,a) \in E_0} \pi(\sigma, a), \tag{2}$$

avec les notations du début de l'introduction et où $E_0$ est un ensemble de couples formés d'une représentation cuspidale unitaire $\sigma$ et d'un entier positif $a$.

Soit $c \in \mathbb{N}$ et soit $\rho_0$ une représentation cuspidale irréductible et unitaire de $\mathrm{GL}(c, F)$. On dit que $\rho_0$ est symplectique (resp. orthogonale) si $L_{v_0}(\rho_0, \wedge^2, s)$ (resp. $L_{v_0}(\rho_0, \mathrm{Sym}^2, s)$) a un pôle en $s = 0$.

On note

$$\mathrm{Jord}_{\rho_0}(\pi_0) := \{a \in \mathbb{N} | \exists (\sigma_0, a) \in E_0; \sigma_0 \simeq \rho_0\} \qquad \text{cf. (2)}$$

cet ensemble étant compté avec multiplicité. On note

- $a_{\rho_0}$ un élément maximal de $\mathrm{Jord}_{\rho_0}$ si cet ensemble est non vide et si cet ensemble est vide,
- $a_{\rho_0} = -1$ si $\rho_0$ et $^L G$ sont tous deux symplectiques ou tous deux orthogonaux,
- $a_{\rho_0} = 0$ si l'un est symplectique et l'autre orthogonal,
- $a_{\rho_0} = \infty$ sinon.

On pose aussi (en abrégant $|\det_{\mathrm{GL}(c,F)}|$ par $|\ |$):

$$\mathrm{Red}_{\rho_0}(\pi_0) := \{s \in \mathbb{R}_{\geq 0} | \rho_0 |\ |^s \times \pi_0 \quad \text{est réductible}\}.$$

On démontre alors:

THÉORÈME. *Soit $c$ et $\rho_0$ comme ci-dessus; alors $\mathrm{Red}_{\rho_0}(\pi_0)$ est vide si et seulement si $a_{\rho_0} = \infty$ et ceci est encore équivalent à ce que $\rho_0$ ne soit pas autoduale. Sinon*

$$\mathrm{Red}_{\rho_0}(\pi_0) = \{(a_{\rho_0} + 1)/2\}.$$

Ce théorème, directement inspiré des travaux de Shahidi, répond aux conjectures que j'avais émises il y a quelques temps.

PROPOSITION. *Pour tout $\rho_0$ comme ci-dessus, l'ensemble $\mathrm{Jord}_{\rho_0}(\pi_0)$ est sans multiplicité et s'il contient un élément $a > 2$ il contient aussi $a - 2$; c'est-à-dire qu'il est sans trou. De plus $\mathrm{Jord}_{\rho_0}(\pi_0)$ a la condition de parité suivante (s'il est non vide): il est formé d'entiers impairs si $\rho_0$ et $^L G$ sont soit tous deux symplectiques soit tous deux orthogonaux et d'entiers pairs sinon. Cet ensemble est donc uniquement déterminé par $\pi_0$ et la représentation $\tilde{\pi}_{v_0}$ ne dépend que de $\pi_0$.*

Expliquons aussi plus précisément ici la normalisation des opérateurs d'entrelacement en ne faisant que les hypothèses (A) et (B); on reprend les notations $\tilde{\pi}_0^L$ et $\tilde{\pi}_0^A$ expliquées au début de l'introduction. Pour $\rho_0$ fixé et $s \in \mathbb{R}$, on note $M(\rho_0 \times \pi_0, s)$ l'opérateur d'entrelacement non normalisé:

$$\rho_0|\ |^s \times \pi_0 \to \rho_0^*|\ |^{-s} \times \pi_0.$$

On définit de façon identique $M(\rho_0^* \times \pi_0, -s)$.

On note $r$ la représentation $\mathrm{Sym}^2\mathbb{C}^c$ si $^LG$ est symplectique et $\wedge^2\mathbb{C}^c$ si $^LG$ est orthogonal.

THÉORÈME.    *Les fonctions méromorphes:*

$$M(\rho_0^* \times \pi_0, -s) \circ M(\rho_0 \times \pi_0, s)$$

*et le produit des facteurs*

$$L_{v_0}(\rho_0 \times \tilde{\pi}_0^A, s) L_{v_0}(\rho_0, r, 2s) L_{v_0}(\rho_0 \times \tilde{\pi}_0^A, 1+s)^{-1} L_{v_0}(\rho_0, r, 1+2s)^{-1}$$
$$\times L_{v_0}(\rho_0^* \times \tilde{\pi}_0^A, -s) L_{v_0}(\rho_0^*, r, -2s) L_{v_0}(\rho_0^* \times \tilde{\pi}_0^A, 1-s)^{-1} L_{v_0}(\rho_0^*, r, 1-2s)^{-1}$$

*diffèrent par une fonction holomorphe inversible. L'opérateur d'entrelacement défini par:*

$$N(\rho_0 \times \pi_0, s) := \bigg( L_{v_0}(\rho_0 \times \tilde{\pi}_0^A, s) L_{v_0}(\rho_0, r, 2s) L_{v_0}(\rho_0 \times \tilde{\pi}_0^A, 1+s)^{-1}$$
$$\cdot L_{v_0}(\rho_0, r, 1+2s)^{-1} \bigg)^{-1} M(\rho_0 \times \pi_0, s)$$

*vérifie la formule du produit sous la forme suivante:*

$$N(\rho_0^* \times \pi_0, -s) N(\rho_0 \times \pi_0, s) \quad \textit{est une fonction holomorphe inversible de } s.$$

*Cet opérateur normalisé $N_{v_0}(\rho_0 \times \pi_0, s)$ est holomorphe pour $s \in \mathbb{C}$ avec $\mathrm{Re}\, s \in \mathbb{R}_{\geq 0}$ mais il peut être identiquement nul. L'ensemble de ses zéros réels positifs ou nuls se calcule explicitement et exclut évidemment $0$. Dans le cas où $\tilde{\pi}_0^A$ est de la forme (2) ils constituent exactement l'ensemble $\{(a-1)/2\}$ pour $a \in \mathrm{Jord}_{\rho_0}(\pi_0)$ avec $a > 1$; ces zéros sont alors simples (il faut (2) pour obtenir la simplicité).*

La normalisation que nous avons définie est donc une normalisation à la Arthur. La normalisation à la Langlands s'obtient (cf. §5) en remplaçant dans les définitions la représentation $\tilde{\pi}_0^A$ par $\tilde{\pi}_0^L$; on démontre alors en plus que l'opérateur d'entrelacement normalisé n'a pas de zéro de partie réelle $\geq 0$. On pourrait d'ailleurs obtenir le même théorème (à la localisation des zéros près) en remplaçant $\tilde{\pi}_0^A$ par $\tilde{\pi}_{v_0}$; c'est une normalisation intermédiaire qui se comprend en considérant que $\tilde{\pi}_{v_0}$ qui dépend du choix de $\tilde{\pi}$ définit un paquet d'Arthur pour $\pi_0$; nous n'avons pas développé ce point de vue dans la mesure où dans l'article on a évité de parler de paquet d'Arthur.

Signalons que Y. Zhang (cf. [Zh]) a retrouvé les résultats concernant les points de réductibilité des induites de cuspidales en prenant d'autres hypothèses que les nôtres; pour sa part elle suppose que ce sont les conjectures de Shahidi (cf. [S$_2$]) qui sont démontrées.

*Remarque.* Supposons que $\pi$ est une représentation cuspidale irréductible de $G(n, \mathbb{A})$ et que $\tilde{\pi}$ est une représentation $G(n)$-discrète de $\mathrm{GL}(2n+\varepsilon, \mathbb{A})$ telle que $\tilde{\pi}_v$ et $\pi_v$ se correspondent par la fonctorialité de Langlands pour tout $v$ excepté un nombre fini de places arbitrairement grand. Alors $\tilde{\pi}$ est autoduale.

Bien sûr cette remarque n'est pas nouvelle; on devrait sans doute même démontrer plus à savoir que $L_v(\tilde{\pi}_v, r, s)$ ($r$ est défini dans l'introduction) a un pôle en $s = 0$ à toute place $v$ mais les arguments ci-dessous sont insuffisants pour cela. On reprend les arguments de [J-S].

Fixons $S$ un nombre fini de places de $k$ telles que $\tilde{\pi}_v$ et $\pi_v$ se correspondent hors de $S$. Alors pour tout $v$ non dans $S$, $\tilde{\pi}_v$ est autoduale. Ainsi $\tilde{\pi}_v^* \simeq \tilde{\pi}_v$ pour tout $v$ non dans $S$. On veut en déduire que $\tilde{\pi}^* \simeq \tilde{\pi}$. C'est bien connu: le support cuspidal de $\tilde{\pi}$ est bien défini et en outre $\tilde{\pi}$ est l'unique représentation $G(n)$-discrète ayant ce support cuspidal. En conjugant dans la chambre de Weyl positive, on écrit uniquement ce support comme un couple $(M, \Sigma)$ où $M$ est un Levi et $\Sigma$ une représentation cuspidale de ce Levi. D'autre part $\tilde{\pi}^*$ est elle aussi de la forme:

$$\times_{(\sigma,a)\in E} \pi(\sigma^*, a).$$

Le support cuspidal de $\tilde{\pi}^*$ est lui aussi bien déterminé et $\tilde{\pi}^*$ est l'unique représentation ayant ce support cuspidal et s'écrivant comme ci-dessus (en autorisant *a priori* un changement de $E$). Le support cuspidal de $\tilde{\pi}^*$ est de la forme $(M, \Sigma^*)$ avec les mêmes notations que ci-dessus (on a encore conjugué dans la chambre de Weyl positive). Il suffit maintenant de démontrer que $\Sigma^* \simeq \Sigma$. On écrit

$$\Sigma =: \times_{(\sigma,x)\in K} \sigma | \, |^x,$$

où $K$ est formé de couples $(\sigma, x)$ avec $\sigma$ cuspidal unitaire irréductible et $x$ demi-entier. On remarque que si $(\sigma, x) \in K$ alors $(\sigma, -x) \in K$ avec la même multiplicité. On a donc aussi

$$\Sigma^* \simeq \times_{(\sigma,x)\in K} \sigma^* | \, |^x.$$

Par induction décroissante sur $x$ on démontre que si $(\sigma, x) \in K$ alors $(\sigma^*, x)$ est aussi dans $K$ avec la même multiplicité que $(\sigma, x)$. Soit $\rho$ une représentation cuspidale unitaire irréductible d'un groupe linéaire:
$$L^S(\rho \times \tilde{\pi}, s) = \times_{(\sigma,x) \in K} L^S(\rho \times \sigma, s-x),$$
$$L^S(\rho \times \tilde{\pi}^*, s) = \times_{(\sigma,x) \in K} L^S(\rho \times \sigma, s-x).$$

En outre ces fonctions $L^S$ ne dépendent que des places hors de $S$ et elles sont donc égales. Si l'on a démontré l'assertion annoncée pour tout $x \geq x_0$ demi-entier fixé, on a:
$$\times_{(\sigma,x) \in K; x < x_0} L^S(\rho \times \sigma, s-x) = \times_{(\sigma,x) \in K; x < x_0} L^S(\rho \times \sigma^*, s-x).$$

On calcule l'ordre du membre de gauche en $s_0 := x_0 + 1/2$ en utilisant des résultats de Jacquet, Shalika [J-S] et de Shahidi [$S_1$] qui assurent que ces fonctions $L^S(\rho \times \sigma, s-x)$ n'ont pas de zéros pour $\operatorname{Re} s - x \geq 1$ et ont un pôle exactement simple en $s = x + 1$ si $\rho \simeq \sigma$. L'ordre cherché vaut donc précisément le cardinal:
$$\{(\sigma, x) \in K | \sigma \simeq \rho, x = x_0 - 1/2\}.$$

D'après ce qui précède:
$$|\{(\sigma, x) \in K | \rho \simeq \sigma, x = x_0 - 1/2\}| = |\{(\sigma, x) \in K | \rho \simeq \sigma^*, x = x_0 - 1/2\}|$$

On obtient alors aisément l'assertion en remplaçant $x_0$ par $x_0 - 1/2$. On montre ainsi que $\Sigma = \Sigma^*$ et d'après ce que l'on a vu au début, $\tilde{\pi} \simeq \tilde{\pi}^*$. Cela termine la preuve.

## 1. Comment neutraliser les places de $k$ différentes de $v_0$

Ici on pose $m := 2n + \varepsilon$ et on fixe une représentation cuspidale unitaire irréductible $\tilde{\pi}$ de $\operatorname{GL}(m, \mathbb{A})$; on fixe aussi $c \in \mathbb{N}$ et $\rho$ une représentation cuspidale unitaire irréductible de $\operatorname{GL}(c, \mathbb{A})$. Pour $v$ une place de $k$, on sait définir [$S_1$] les facteurs $L_v(\rho_v \times \tilde{\pi}_v, s)$ comme fonction méromorphe de $s$ ainsi que les facteurs $L_v(\rho_v, r, s)$ où $r$ a été défini dans l'introduction. On fixe $S$ un ensemble fini de places de $k$ contenant les places archimédiennes et arbitrairement grand et contenant $v_0$. On veut le lemme suivant:

LEMME 1. *Quitte à augmenter $S$, il existe un caractère unitaire $\chi$ de $k^* \backslash \mathbb{A}^*$, non ramifié hors de $S$ et tel que:*
- $\chi_{v_0} = 1$ *pour tout $v \in S$, les facteurs $L_v(\rho_v \otimes \chi_v \times \tilde{\pi}_v, s)$ et $L_v(\rho_v \otimes \chi_v, r, s)$ sont holomorphes non nuls pour tout $s \in \mathbb{R}$,*
- *les fonctions $L(\rho \otimes \chi \times \tilde{\pi}, s)$ et $L(\rho \otimes \chi, r, s)$ sont holomorphes pour tout $s \in \mathbb{R}$.*

On fixe *a priori* $S$ contenant $v_0$, les places archimédiennes et tel que $\rho$ et $\tilde{\pi}$ sont non ramifiées hors de $S$. On note $p_0$ le nombre premier divisant le nombre d'éléments du corps résiduel en la place $v_0$ et on note $S'$ l'ensemble des places non archimédiennes de $k$ pour lesquelles le nombre d'éléments du corps résiduel est multiple de $p_0$. On suppose que $S$ contient $S'$.

On va commencer par faire un choix de $\chi$ qui évite les ennuis aux places de $S' - \{v_0\}$. Pour cela on fixe une telle place $v$; les pôles sur l'axe réel de $L_v(\rho_v \times \tilde{\pi}_v, s)$ viennent des entrelacements entre le support cuspidal de $\rho_v||^s$ et celui de $\tilde{\pi}_v \simeq \tilde{\pi}_v^*$; on peut éviter ces entrelacements en tensorisant $\rho_v$ par un caractère $\chi_v$ d'ordre fini. En fixant convenablement $\chi_v$ on évite aussi tout entrelacement entre le support cuspidal de $\rho_v||^s \otimes \chi_v$ et $\rho_v^*||^{-s} \otimes \chi_v^{-1}$ ce qui évite d'avoir des pôles aux facteurs $L_v(\rho_v \otimes \chi_v, r, s)$ sur l'axe réel.

On note $\chi_{S'}$ le produit des choix $\chi_v$ faits pour $v \in S' - v_0$ et du caractère trivial en $v_0$. On obtient ainsi un caractère d'ordre fini de $\times_{v \in S'} k_v$. On a besoin du lemme suivant, dont la démonstration m'a été fournie par Waldspurger:

LEMME 2.   *Le caractère $\chi_{S'}$ se prolonge en un caractère lisse de $\mathbb{A}^*$ trivial sur $k^*$.*

On note $N$ un entier tel que $\chi_{S'}$ soit d'ordre divisant $N$. Soit ici $S$ un ensemble fini de places de $k$ contenant $S'$ ainsi que les places archimédiennes et celles dans lesquelles $N$ se ramifie. Pour tout $v$ dans $S - S'$ on fixe un caractère lisse de $k_v^*$ d'ordre divisant $N$, noté $\chi_v$ par exemple le caractère trivial. On construit ci-dessous un caractère lisse $\chi$ des idèles $I_k$ de $k$ trivial sur $k^* I_k^N$ et vérifiant pour tout $v \in S$:

$$\chi_{|k_v^*} = \chi_v. \tag{1}$$

Bien sûr on ne sait rien sur la ramification de $\chi$ hors de $S$. On fixe $S_1$ un ensemble fini de places de $k$ disjoint de $S$. On remarque que le sous-groupe de $I_k$:

$$H_{S_1} := \big(\prod_{v \in S} k_v^* \prod_{v \notin S \cup S_1} O_v^*\big) I_k^N$$

contient un sous-groupe compact ouvert de $I_k$; on utilise ici le fait que pour $v \in S_1$, $1 + \omega_v O_v$ est inclus dans $k_v^N$ car $N$ ne se ramifie pas à cette place.

Il existe un caractère $\chi'$ de $H_{S_1}$ vérifiant (1) ci-dessus et trivial sur $(\prod_{v \notin S \cup S_1} O_v^*) I_k^N$. Il faut maintenant montrer que quitte à augmenter $S_1$ (ce qui diminue $H_{S_1}$), on peut prolonger $\chi'$ en un caractère de $H_{S_1} k^*$ trivial sur $k^*$.

On sait que $E := k^* \cap \big(\prod_{v \in S} k_v^* \prod_{v \notin S} O_v^*\big)$ est de type fini. On peut donc appliquer le théorème 1 de [C]: il existe un ensemble fini de places $S_1$ disjoint de $S$ et pour tout $v \in S_1$ un compact ouvert $U_v$ de $O_v^*$ tel que:

$$k^* \cap \big(\prod_{v \in S} k_v^* \prod_{v \in S_1} U_v \prod_{v \notin S \cap S_1} O_v^*\big) \subset k^{*N}.$$

On vérifie alors que ce choix de $S_1$ convient. Comme le caractère est à valeurs dans un corps algébriquement clos, il n'y a aucun problème pour le prolonger de $k^*H_{S_1}$ à $I_k$. Cela termine la preuve du lemme.

On poursuit la preuve du lemme 1. Ici on fixe $\chi_1$ un caractère adélique de $k$ ayant la restriction souhaitée à chaque place de $S'$. Maintenant on augmente $S$ de façon à y inclure les places où $\chi_1$ est ramifié. Et on va démontrer que l'on peut construire $\chi_2$ un caractère non ramifié des idèles de $k$ tel que le lemme soit satisfait pour $\chi = \chi_1\chi_2$. On construit $\chi_2$ sous la forme $||^{iz/\log p_0}$ où $p_0$ a été défini au début de la preuve; ces caractères sont partout non ramifiés et sont triviaux en toutes places de $S'$. Après tensorisation par $\chi_2$ on a évité les entrelacements, pour tout $s \in \mathbb{R}$, entre $\rho||^s \otimes \chi_1\chi_2$ et $\tilde{\pi}$ ou $\rho^*||^{-s} \otimes \chi_1^{-1}\chi_2^{-1}$ pour tout choix de $z$ sauf un nombre fini; ce sont ces entrelacements qui produisent les pôles des fonctions $L$ globales. Par le même procédé, on élimine les ennuis éventuels aux places archimédiennes et aux places dans $S - S'$. Cela termine la preuve.

## 2. Calcul de $\mathrm{Red}_{\rho_0}(\pi_0)$ dans le cas le plus simple

2.1. *Quelques généralités.* Dans toute la suite de cet article, on fixe $\pi_0$ comme dans l'introduction ainsi que $\tilde{\pi}$ une représentation irréductible de $\mathrm{GL}(2n + \varepsilon, \mathbb{A})$ et $\pi$ une représentation cuspidale irréductible de $G(n)$ satisfaisant les hypothèses (A), (B), (C) de l'introduction; on suppose que cela est possible. On pose

(1) $$\tilde{\pi}_0 := \tilde{\pi}_{v_0}.$$

On fixe aussi $c \in \mathbb{N}$ et $\rho_0$ une représentation cuspidale unitaire irréductible de $\mathrm{GL}(c, F)$.

On construit (cf. [He, App.]) une représentation cuspidale, unitaire irréductible, $\rho$, de $\mathrm{GL}(c, \mathbb{A})$. On fixe $S$ un ensemble fini de places contenant les places archimédiennes et tel que hors de $S$, $\rho, \pi, \tilde{\pi}$ sont non ramifiés et $\tilde{\pi}$ et $\pi$ se correspondent par la fonctorialité de Langlands. Quitte à augmenter $S$ on fixe un caractère $\chi$ (cf. §1) tel que $\chi_{v_0} = 1$, $\chi$ est non ramifié hors de $S$ et les fonctions $L(\rho \otimes \chi \times \tilde{\pi}, s)$ et $L(\rho \otimes \chi, r, s)$ sont holomorphes en tout point $s \in \mathbb{R}$ et les fonctions $L_v(\rho \otimes \chi \times \tilde{\pi}, s)$ et $L_v(\rho \otimes \chi, r, s)$ sont holomorphes non nulles en tout $s \in \mathbb{R}$. Evidemment pour simplifier les notations, on va remplacer $\rho$ par $\rho \otimes \chi$.

Soit $f$ dans l'induite automorphe $\rho \times \pi$. On définit l'opérateur d'entrelacement non normalisé:

$$M(\rho \times \pi, s)f : \rho||^s \times \pi \to \rho^*||^{-s} \times \pi,$$

en tant que fonction méromorphe de $s$; cela peut aussi se faire localement. On exclut ici (puisque l'on travaille avec les composantes connexes, cf. (A)) le cas où $G = O(2n, F)$ avec $n = 0$ et $c$ impair qui sera vu avant la proposition de 2.2. On suppose que $f$ est totalement décomposée et est non ramifiée hors de $S$. On a une égalité de fonctions méromorphes:

$$M(\rho \times \pi, s)f = r^S(\rho \times \pi, s)M_S(\rho \times \pi, s)f_S,$$

où l'exposant veut dire que l'on évite les places indiquées, l'indice veut dire que l'on fait le produit sur les places indiquées et où (cf. [L])

$$r^S(\rho \times \pi, s) = L^S(\rho \times \pi, s)L^S(\rho, r, 2s)L^S(\rho \times \pi, s+1)^{-1}L^S(\rho, r, 2s+1)^{-1}.$$

On fait des constructions analogues en fixant $f^*$ dans l'induite automorphe de $\rho^* \times \pi$.

En tordant $\rho$ par un caractère comme dans Section 1 on impose en plus que pour tout $v \in S - \{v_0\}$, $M_v(\rho \times \pi, s)f$ et $M_v(\rho^* \times \pi, s)f^*$ est holomorphe pour tout $s \in \mathbb{R}$; ici on doit éviter les entrelacements entre le support cuspidal de $\rho_v||^s$ et celui de $\pi_v$. J'ai conscience que ma présentation risque d'induire en confusion. On procède de la façon suivante, avec la notation $S'$ de la preuve du lemme 1: on fixe $\chi_{S'}$ caractère de $\times_{v \in S'} k_v$ de façon à éviter les pôles pour les facteurs $L_v$ et pour les opérateurs d'entrelacement ci-dessus pour tout $v \in S' - \{v_0\}$. Puis on étend $\chi_{S'}$ en un caractère adélique; on fixe alors $S$ de façon à ce que tout soit non ramifié hors de $S$ et on construit $\chi_2$ un caractère partout non ramifié trivial sur $k_v^*$ pour tout $v \in S'$ (de la forme $||^{iz/\log p_0}$, qui évite:

les pôles des facteurs $L_v$ et des opérateurs d'entrelacement aux places de $S - S'$,
les pôles des fonctions $L$ globales,
et tel que $\rho \otimes \chi_1 \otimes \chi_2 \not\simeq \rho^* \otimes \chi_1^{-1} \otimes \chi_2^{-1}$.

On remplace alors $\rho$ par $\rho \otimes \chi_1 \otimes \chi_2$; on n'a rien changé à la place $v_0$. La dernière propriété écrite ci-dessus assure que $\rho \not\simeq \rho^*$.

*Remarque.* Avec les hypothèses ci-dessus, il n'existe pas de forme automorphe de carré intégrable dont le support cuspidale soit de la forme $\rho||^{s_0} \times \pi$ avec $s_0 \in \mathbb{R}$ et qui soit non ramifiée hors de $S$.

Cela est du à Langlands et est bien connu dans le cas où le support cuspidal est un parabolique maximal et une représentation cuspidale d'un Levi de ce parabolique. Pour l'expliquer, on va introduire des notations qui nous serviront dans la suite.

On note $P$ le sous-groupe parabolique maximal de $G(n+c)$ de Levi $GL(c, \rho) \times G(n)$; si $n = 0$ et $G$ est un groupe orthogonal paire, il y a une petite différence, $P$ n'étant pas uniquement défini (à conjugaison près) par son Levi mais on laisse ce cas de côté pour le reprendre plus loin à la fin de 2.2.

Pour $f$ et $f^*$ comme au début de ce paragraphe:
$$E(f,s)_P = f_s + r^S(\rho \times \pi, s) \, M_S(\rho \times \pi, s) f_s,$$
$$E(f^*, -s)_P = f^*_{-s} + r^S(\rho^* \times \pi, -s) \, M_S(\rho^* \times \pi, -s) f^*_{-s}$$

et ce sont les seuls termes constants non nuls donnant le support cuspidal de ces séries d'Eisenstein.

Langlands a démontré qu'il ne peut y avoir de formes automorphes de carré intégrable avec comme support cuspidal $\rho| \, |^{s_0} \times \pi$ ($s_0 \in \mathbb{R}_{\geq 0}$, le support n'est défini qu'à conjugaison près) et non ramifiée hors de $S$ que si les séries d'Eisenstein $E(f,s)$ construite avec $f$ non ramifiée hors de $S$ ont un pôle en $s = s_0$. Et alors (cf. [M-W, 5.3.8])

$$w_0(\rho| \, |^{s_0} \otimes \pi) \simeq \rho| \, |^{-s_0} \otimes \pi,$$

où $w_0$ est l'élément non trivial du groupe de Weyl relatif à l'action du centre du Levi de P. Cette relation est encore équivalente à $\rho \simeq \rho^*$ ce que nous avons exclu.

2.2. *Le cas simple de l'hypothèse* (C). Dans ce paragraphe, on suppose que $\pi_0$ et $\tilde{\pi}$ fixés ci-dessus vérifient aussi l'hypothèse (C) en plus de (A) et (B). Donc par hypothèse, il existe un ensemble $E_0$ formé de couple $(\sigma, a)$ où $\sigma$ est une représentation cuspidale unitaire d'un groupe linéaire et $a \in \mathbb{N}$ tel que:

(2) $$\tilde{\pi}_0 \simeq \times_{(\sigma, a) \in E_0} \pi(\sigma, a).$$

On rappelle que $\pi(\sigma, a)$ est l'unique quotient irréductible de l'induite:
$$\times_{j \in [(a-1)/2, -(a-1)/2]} \sigma |\det|^j.$$

Soit $s \in \mathbb{R}$; on note: $[\rho_0| \, |^s : \tilde{\pi}_0]$ la multiplicité de $\rho_0| \, |^{s_0}$ dans le support cuspidal de $\tilde{\pi}_{v_0}$; en termes concrets et en reprenant la notation (2):
$$[\rho_0| \, |^s : \tilde{\pi}_0] := |\{(\sigma, a) \in E_0 | \sigma \simeq \rho_0, s \in \{-(a-1)/2 + k\}_{k \in [0, a-1]}\}|.$$

Cette définition peut se faire en remplaçant $\rho$ par $\rho^*$ et remarquons que l'autodualité de $\tilde{\pi}_0$ entraîne que:
$$[\rho_0^*| \, |^{-s} : \tilde{\pi}_0] = [\rho_0| \, |^s : \tilde{\pi}_0].$$

LEMME. *Avec les hypothèses ci-dessus*: *soit $s_0 \in \mathbb{R}$; alors*:
$$\operatorname{ordre}_{s=s_0} r^S(\rho \times \pi, s) = [\rho_0| \, |^{s_0} : \tilde{\pi}_0] - [\rho_0| \, |^{s_0+1} : \tilde{\pi}_0] + \varepsilon(s_0) + \eta(s_0)$$

*où $\varepsilon(s_0) := \operatorname{ordre}_{s=s_0} L_{v_0}(\rho_0, r, 2s)^{-1} L_{v_0}(\rho_0, r, 2s+1)$ et où*
$$\eta(s_0) := \operatorname{ordre}_{s=s_0} L(\rho \times \tilde{\pi}, s) L(\rho, r, 2s) L(\rho \times \tilde{\pi}, s+1)^{-1} L(\rho, r, 2s+1)^{-1};$$

*et avec les mêmes notations*
$$\operatorname{ordre}_{s=s_0} r^S(\rho^* \times \pi, -s) = [\rho| \, |^{s_0} : \tilde{\pi}_0] - [\rho| \, |^{s_0-1} : \tilde{\pi}_0] + \varepsilon(-s_0) - \eta(s_0).$$

On a:
$$L^S(\rho \times \tilde{\pi}, s) = \prod_{v \in S} L_v(\rho_v \times \tilde{\pi}_v, s)^{-1} L(\rho \times \tilde{\pi}, s).$$

D'où, avec les hypothèses faites:

$$\text{ordre}_{s=s_0} L^S(\rho \times \tilde{\pi}, s) = \text{ordre}_{s=s_0} L_{v_0}(\rho_0 \times \tilde{\pi}_0, s)^{-1} + \text{ordre}_{s=s_0} L(\rho \times \tilde{\pi}, s).$$

Or:
$$\text{ordre}_{s=s_0} L_{v_0}(\rho_0 \times \tilde{\pi}_0, s)^{-1} = [\rho_0 |\ |^{s_0} : \tilde{\pi}_0],$$

d'où le premier calcul de l'énoncé. Le deuxième est du même type: l'autodualité de $\tilde{\pi}_v$ en toute place assure que si $L_v(\rho_v \times \tilde{\pi}_v, s)$ est holomorphe pour tout $s \in \mathbb{R}$, il en est de même de $L_v(\rho_v^* \times \tilde{\pi}_v, s)$. Et pour trouver l'analogue de $\eta(s_0)$ on utilise l'équation fonctionnelle pour avoir l'égalité:

$$\eta(s_0) = \text{ordre}_{s=s_0} L(\rho^* \times \tilde{\pi}, 1-s) L(\rho^*, r, 1-2s) L(\rho^* \times \tilde{\pi}, -s)^{-1} L(\rho^*, r, -2s)^{-1}.$$

Cela termine la preuve.

On reprend les notations $\text{Jord}_{\rho_0}(\pi_0)$ de l'introduction. Pour $s_0 \in \mathbb{R}_{\geq 0}$, on pose:
$$n_1(s_0) := |\{a \in \text{Jord}_{\rho_0}(\pi_0) | s_0 = (a+1)/2\}| + \varepsilon_0,$$

où $\varepsilon_0$ vaut 0 pour tout $s_0$ sauf éventuellement pour $s_0 = 1/2$ où $\varepsilon_0$ vaut 1 exactement quand $L_{v_0}(\rho, r, s)$ a un pôle en $s = 0$. Ainsi $n_1(1/2) = \varepsilon_0$. On pose aussi pour $s_0 \in \mathbb{R}_{\geq 0}$:
$$n_0(s_0) := |\{a \in \text{Jord}_{\rho_0}(\tilde{\pi}_0) | s_0 = (a-1)/2\}|.$$

En particulier pour $s_0$ non demi-entier ces nombres sont nuls et si $s_0 \geq 1$:
$$n_1(s_0) = n_0(s_0 - 1).$$

Pour $s_0 \in\ ]0, 1[$, $n_1(s_0) \leq 1$.

*Remarque.* Pour tout $s_0 \in \mathbb{R}, s_0 \geq 0$:
$$[\rho |\ |^{s_0} : \tilde{\pi}_0] - [\rho |\ |^{s_0+1} : \tilde{\pi}_0] = n_0(s_0)$$

et pour $s_0 > 0$
$$[\rho |\ |^{s_0} : \tilde{\pi}_0] - [\rho |\ |^{s_0-1} : \tilde{\pi}_0] = -n_1(s_0) + \varepsilon_0.$$

Quant à $\eta(s_0)$ défini dans le lemme ci-dessus, c'est un nombre positif ou nul pour tout $s_0 > 0$ et nul pour $s_0 = 0$.

De plus $\varepsilon(s_0) = 0$ pour tout $s_0 \in \mathbb{R} - \{0\}$ sauf éventuellement en $s_0 = -1/2$ où $\varepsilon(-1/2) = \varepsilon_0$.

Les deux premières assertions sont un calcul explicite en utilisant la description du support cuspidal de $\tilde{\pi}_0$. La troisième assertion vient des propriétés des facteurs $L_{v_0}(\rho_0, r, s)$: holomorphie sans zéro pour tout $s \in \mathbb{R} - \{0\}$ et pôle au plus simple en $s = 0$.

La propriété de $\eta(s_0)$ utilise le calcul explicite suivant qui reprend les notations (1) de l'introduction:

$$L(\rho \times \tilde{\pi}, s)/L(\rho \times \tilde{\pi}, s+1)$$
$$= \prod_{(\sigma,a)\in E} \prod_{k\in[0,a-1]} L(\rho \times \sigma, s + (a-1)/2 - k)/L(\rho \times \sigma, s + (a-1)/2 - k + 1)$$
$$= \prod_{(\sigma,a)\in E} L(\rho \times \sigma, s - (a-1)/2)/L(\rho \times \sigma, s + (a+1)/2).$$

Or les numérateurs et dénominateurs ci-dessus n'ont pas de pôles par la propriété de $\rho$ (ils créeraient des pôles pour $L(\rho \times \tilde{\pi}, s)$) par contre ils peuvent avoir des zéros. Mais si $\operatorname{Re} s \geq 0$, les dénominateurs n'ont pas de zéros puisque les zéros pour les fonctions $L$ de paires de représentations cuspidales unitaires sont dans la bande critique $[0,1]$. De même $L(\rho, r, 2s)/L(\rho, r, 2s+1)$ est une fonction holomorphe pour $s \in \mathbb{R}_{\geq 0}$ avec éventuellement des zéros pour $s \in [0, 1/2]$. En tout cas on a la positivité pour $\eta(s_0)$. La dernière propriété résulte des définitions.

On définit $\operatorname{irr}(s_0) = 0$ si $s_0 \notin \operatorname{Red}_{\rho_0}(\pi)$ et $\operatorname{irr}(s_0) = 1$ sinon. Evidemment cette définition dépend de $\rho_0$ et de $\pi_0$; mais comme c'est une notation intermédiaire, on n'a pas marqué cette dépendance dans la notation. On pose $\varepsilon'_0 = 0$ si $L_v(\rho_0, r, s)$ n'a pas de pôle en $s = 0$ et 1 dans le cas contraire.

Le cas de $G = \operatorname{O}(2n, F)$ pose un problème spécifique que l'on va aborder ici, le but étant de montrer dans quel cas on peut confondre $\operatorname{O}(2n, F)$ et $\operatorname{SO}(2n, F)$. On va montrer que c'est le cas soit si $c$ est pair soit si $c$ est impair mais la restriction de $\pi_0$ à $\operatorname{SO}(2n, F)$ est irréductible.
Remarquons que le parabolique $P$ est conjugué du parabolique opposé $\overline{P}$ dans $\operatorname{SO}(2n+2c, F)$ si et seulement si $n \neq 0$ ou $n = 0$ et $c$ est pair. Supposons $n \neq 0$: si la restriction de $\pi_0$ à $\operatorname{SO}(2n, F)$ est irréductible alors $\rho_0||^{s_0} \times \pi_0$ en tant qu'induite à $\operatorname{O}(2n+2c, F)$ est irréductible si et seulement si ceci est vrai en tant qu'induite à $\operatorname{SO}(2n+2c, F)$. Dans ce cas, il n'y a pas d'ambiguité à confondre $\pi_0$ en une représentation de $\operatorname{SO}(2n, F)$. Par contre si la restriction de $\pi_0$ à $\operatorname{SO}(2n, F)$ n'est pas irréductible fixons $\pi'_0$ une composante irréductible, alors l'induite de $\rho_0||^{s_0} \times \pi'_0$ à $\operatorname{O}(2n+2c, F)$ est naturellement isomorphe à l'induite de $\rho_0||^{s_0} \times \pi_0$ à $\operatorname{O}(2n+2c, F)$. Si $c$ est pair, l'induite de $\rho_0||^{s_0} \times \pi'_0$ à $\operatorname{SO}(2n+2c, F)$ à $\operatorname{O}(2n+2c, F)$ est irréductible si et seulement si cela est vrai pour l'induite à $\operatorname{SO}(2n+2c, F)$ (on le vérifie sur les modules de Jacquet). Dans ce cas aussi, on peut donc se ramener sans grandes précautions à $\operatorname{SO}(2n, F)$. Par contre si $c$ est impair la situation est différente, il n'y a pas à invoquer de gros théorème pour savoir que l'induite de $\rho_0||^{s_0} \times \pi'_0$ à $\operatorname{SO}(2n+2c, F)$ est irréductible pour tout $s_0 \in \mathbb{R}$. Quant à l'induite à $\mathrm{O}(2n+2c,F)$ elle est irréductible si et seulement si $\rho_0||^{s_0} \not\simeq \rho_0^*||^{-s_0}$. D'où, sous l'hypothèse faite: $G = \mathrm{O}(2n,F)$, $c$ impair et la restriction de $\pi_0$ à $\mathrm{SO}(2n,F)$ non irréductible (en particulier $n \neq 0$) un calcul immédiat donne:

$$\mathrm{irr}(s_0) = 0 \ \forall \ s_0 \in \mathbb{R}_{>0},$$
$$\mathrm{irr}(0) = 1 \Leftrightarrow \rho_0 \simeq \rho_0^*.$$

Reste à voir le cas où $n = 0$: si $c$ est pair, $P$ est conjugué de $\overline{P}$ dans $\mathrm{SO}(2c,F)$ et un élément bien choisi de $\mathrm{O}(2c,F)$ conjugue $P$ en l'autre parabolique maximal, $Q$, de Levi isomorphe à $\mathrm{GL}(c,F)$. Alors l'induite de $\rho_0||^{s_0}$ à $\mathrm{SO}(2c,F)$ n'est pas invariante par l'action de $\mathrm{O}(2c,F)$; plus précisément le module de Jacquet de cette induite relativement à $Q$ est nul alors que le module de Jacquet de son image par un élément de $\mathrm{O}(2c,F)$ conjugant $P$ et $Q$ a la propriété opposée d'avoir un module de Jacquet relativement à $P$ nul. On obtient donc que l'induite à $\mathrm{O}(2c,F)$ est irréductible si et seulement si cela est vrai quand on induit à $\mathrm{SO}(2c,F)$. Par contre si $c$ est impair, on obtient comme ci-dessus:

$$\mathrm{irr}(s_0) = 0 \ \forall \ s_0 \in \mathbb{R}_{>0},$$
$$\mathrm{irr}(0) = 1 \Leftrightarrow \rho_0 \simeq \rho_0^*.$$

La proposition suivante va faire rentrer le groupe orthogonal dans le rang.

PROPOSITION. *Soit $s_0 \in \mathbb{R}, s_0 > 0$. Alors*:

$$n_0(s_0) = n_1(s_0) - \mathrm{irr}(s_0).$$

*En $s_0 = 0$: Si $\rho_0 \simeq \rho_0^*$*:

$$n_0(0) + \mathrm{irr}(0) + \varepsilon_0' = 1,$$

*sinon*:

$$n_0(0) + \mathrm{irr}(0) + \varepsilon_0' = 0.$$

*où $\varepsilon_0'$, rappelons le, $\varepsilon_0' = 0$ si $L_v(\rho_0, r, s)$ n'a pas de pôle en $s = 0$ et $1$ dans le cas contraire.*

On suppose d'abord que l'on n'est pas dans le cas où $n = 0$ et $c$ est impair et encore que $s_0 > 0$. On pose

$$\tilde{n}_0(s_0) := n_0(s_0) + \eta(s_0),$$
$$\tilde{n}_1(s_0) := n_1(s_0) + \eta(s_0).$$

Et

$$\delta(s_0) := \tilde{n}_0(s_0) - \tilde{n}_1(s_0) + \mathrm{irr}(s_0).$$

D'après ce que l'on a vu pour $s_0 \geq 0$ et $s_0 \in \mathbb{R}$:

$$\mathrm{ordre}_{s=s_0} r^S(\rho \times \pi, s) = \tilde{n}_0(s_0).$$

Or
$$\text{ordre}_{s=s_0} r^S(\rho^* \times \pi, -s) = -\tilde{n}_1(s_0)$$

On reprend les notations de 2.1. Si $s_0$ est un point de réductibilité, on fixe $f^*$ tel que $f^*_{-s_0}$ soit dans le noyau de l'opérateur $M_{v_0}(\rho_0^* \times \pi_0, -s_0)$ sinon on prend $f^*$ quelconque et on pose:

$$\phi_s := (s-s_0)^{\tilde{n}_1(s_0)-\text{irr}(s_0)} M(\rho^* \times \pi, -s) f^*_{-s}.$$

C'est une fonction holomorphe de $s$ à valeurs dans les sections du fibré trivial des représentations $\rho| \ |^s \times \pi$ non ramifiées hors de $S$. Comme les zéros de $M_{v_0}(\rho^* \times \pi, -s)$ sont au plus simples d'après les résultats d'Harish-Chandra et que la représentation $\rho| \ |^{s_0} \times \pi_0$ est de longueur au plus 2, on voit que l'ensemble des valeurs $\phi_{s_0}$ obtenues en faisant varier $f^*$ (avec la propriété imposée) décrit toute l'induite $\rho| \ |^{s_0} \times \pi$. L'équation fonctionnelle montre que

$$\left( E(f^*, -s) - (s-s_0)^{-\tilde{n}_1(s_0)+\text{irr}(s_0)} E(\phi_s, s) \right) = 0.$$

On calcule le terme constant suivant $P$ de cette fonction et on trouve:
$$= f^*_{-s} + (s-s_0)^{-\tilde{n}_1(s_0)+\text{irr}(s_0)} \phi_s - (s-s_0)^{-\tilde{n}_1(s_0)+\text{irr}(s_0)} \phi_s$$
$$+ (s-s_0)^{\delta(s_0)} \left( r^S(\rho \times \pi, s)/(s-s_0)^{\tilde{n}_0(s_0)} \right) M_S(\rho \times \pi, s) \phi_s.$$

Pour un bon choix de $f^*$, $\left( r^S(\rho \times \pi, s)/(s-s_0)^{\tilde{n}_0(s_0)} \right) M_S(\rho \times \pi, s) \phi_s$ est non nul (cf. ci-dessus), d'où $\delta(s_0) = 0$ ce qui entraîne la première partie de la proposition.

Situation en $s = 0$. D'après Harish-Chandra on sait que $M_{v_0}(\rho_0 \times \pi_0, s)$ a un pôle en $s = 0$ si et seulement si $\rho_0 \times \pi_0$ est irréductible et est isomorphe à sa conjuguée par l'élément du groupe de Weyl qui conjugue $P$ et $\overline{P}$. C'est-à-dire qu'il faut $\rho_0 \simeq \rho_0^*$ et en plus, si $c$ est impair et $G = \text{O}(2n)$ que la restriction de $\pi_0$ à $\text{SO}(2n, F)$ soit irréductible (on rappelle que l'on a supposé que si $c$ est impair $n \neq 0$). Quand il y a pôle, le pôle est simple. Pour le groupe orthogonal $\text{O}(2n, F)$ ce sont exactement les cas où l'on a le droit de travailler avec $\text{SO}(2n, F)$ plutôt que $\text{O}(2n, F)$.

D'autre part l'opérateur d'entrelacement global $M(\rho \times \pi, s)$ est holomorphe unitaire en $s = 0$. Or pour un bon choix de $f$:

$$0 = \text{ordre}_{s=0} M(\rho \times \pi, s) f = \text{ordre}_{s=0} r^S(\rho_0 \times \pi_0, s) + \text{ordre}_{s=0} M_{v_0}(\rho_0 \times \pi_0, s) f_{v_0}.$$

Supposons d'abord que $\rho_0 \simeq \rho_0^*$ et que si $G = \text{O}(2n)$ soit $c$ est pair soit la restriction de $\pi_0$ à $\text{SO}(2n, F)$ est irréductible avec $n \neq 0$. On obtient donc

$$0 = n_0(0) + \varepsilon'_0 - 1 + \text{irr}(0).$$

Cela force $n_0(0) + \varepsilon_0' + \mathrm{irr}(0) = 1$, ce qui est l'assertion cherchée. Dans les cas restants, on a tout simplement l'irréductibilité de l'induite et l'égalité:

$$n_0(0) + \varepsilon_0' + \mathrm{irr}(0) = 0.$$

Il reste à voir le cas où $G = \mathrm{O}(2n, F)$ avec $c$ impair et $n = 0$. Dans ce cas l'induite à $\mathrm{SO}(2c, F)$ de $\rho_0$ est irréductible par contre l'induite à $\mathrm{O}(2c, F)$ est réductible si et seulement si $\rho_0 \simeq \rho_0^*$. L'argument avec l'opérateur d'entrelacement donne $\varepsilon_0' = 0$ et on obtient donc encore la proposition.

En un point $s_0 \neq 0$, on considère l'opérateur d'entrelacement associé à un élément, $w_0$, du groupe de Weyl non trivial qui laisse stable le Levi $\mathrm{GL}(c, F)$, $M(w_0, s_0)$; il est défini pour le groupe $\mathrm{SO}(2c)$. Il est holomorphiquement défini localement et globalement. Et $M(w_0, -s_0)M(w_0, s_0) = 1$. On applique alors la première partie de la démonstration pour obtenir $\mathrm{irr}(s_0) = 0$ pour tout $s_0$ où l'induite a lieu pour $\mathrm{SO}(2c, F)$. La représentation induite n'est pas stable sous l'action de $\mathrm{O}(2c, F)$ quand $s_0 \neq 0$ d'où aussi l'irréductibilité pour l'induite à $\mathrm{O}(2c, F)$ et *a fortiori* la proposition.

2.3. *Preuve du théorème de l'introduction.* On rappelle l'énoncé: $\pi_0$ et $\rho_0$ sont fixés. On a défini $\mathrm{Red}_{\rho_0}(\pi_0)$ et $\mathrm{Jord}_{\rho_0}(\pi_0)$ ce qui suppose l'existence de $\tilde{\pi}_0$ (cf. l'introduction). On renvoie aussi à l'introduction pour la définition de $a_{\rho_0}$. Alors:
- L'ensemble $\mathrm{Red}_{\rho_0}(\pi_0)$ est non vide si et seulement si $\rho_0 \simeq \rho_0^*$, ce qui est encore équivalent à $a_{\rho_0} \neq \infty$ si $\mathrm{Red}_{\rho_0}(\pi_0)$ est non vide il est précisément réduit au point $(a_{\rho_0} + 1)/2$.
- L'ensemble $\mathrm{Jord}_{\rho_0}(\pi_0)$ est sans multiplicité, sans trou et a la condition de parité (cf. l'introduction et ci-dessous).

Rappelons que la condition de parité est la suivante: soit $L_{v_0}(\rho_0, r, s)$ a un pôle en $s = 0$ est $\mathrm{Jord}_{\rho_0}(\pi_0)$ est formé d'entiers pairs, soit il n'y a pas de pôle (alors que $\rho_0 \simeq \rho_0^*$) et $\mathrm{Jord}_{\rho_0}(\pi_0)$ est formé d'entier impair.

On fixe $\rho_0$ comme ci-dessus et la proposition dit que pour tout $s_0 \in \mathbb{R}_{>0}$:

$$n_0(s_0) + \mathrm{irr}(s_0) = n_1(s_0).$$

Rappelons que par définition pour tout $s_0 \geq 1$, $n_1(s_0) = n_0(s_0 - 1)$. De plus $n_0(s_0)$ est nul pour $s_0$ grand. Ainsi pour tout $s_0 \in \mathbb{R}_{>0}$:

$$\tag{1} \sum_{k \in \mathbb{N}} \mathrm{irr}(s_0 + k) = n_1(s_0).$$

Par définition $n_1(s_0) \leq 1$ pour $s_0 \in ]0, 1[$ avec égalité si et seulement si $s_0 = 1/2$ et $L_{v_0}(\rho_0, r, s)$ a un pôle en $s = 0$. De plus $n_1(1) = n_0(0)$. Et d'après la proposition ci-dessus $n_0(0) \leq 1$ avec égalité si et seulement si 0 n'est pas point de réductibilité et $L_{v_0}(\rho, r, s)$ n'a pas de pôle en $s = 0$ et $\rho_0 \simeq \rho_0^*$.

(2) On suppose que $\text{Jord}_{\rho_0}(\pi_0)$ est non vide. Alors on note $a'$ un élément de cet ensemble tel que $a' + 2$ n'est pas dans l'ensemble. On pose $s_0 := (a'+1)/2$. Alors $n_0(s_0) = 0$ et $n_1(s_0) \geq 1$. La proposition entraîne que $s_0 \in \text{Red}_{\rho_0}(\pi_0)$ et que $n_1(s_0) = 1$.

Si $\rho_0$ n'est pas isomorphe à $\rho_0^*$ alors $\text{irr}(s_0 + k) = 0$ pour tout $s_0 \geq 0$ et $k \geq 0$ d'où avec (1) $n_0(s_0) = 0$; ceci est aussi vrai en $s_0 = 0$ par la proposition ci-dessus. Ainsi $\text{Red}_{\rho_0}(\pi_0) = \emptyset$ et $\text{Jord}_{\rho_0}(\pi_0) = \emptyset$.

Supposons maintenant que $\rho_0$ soit autodual et en plus que $L_{v_0}(\rho_0, r, s)$ a un pôle en $s = 0$. Si $G = \text{O}(2n, F)$ cela force $c$ à être pair (cf. la proposition ci-dessus). Sous cette hypothèse $n_0(0) = 0$ et $0$ n'est pas un point de réductibilité. On a $n_1(1/2) = 1$ et $n_1(x) = 0$ pour tout $x \in ]0,1[-\{1/2\}$. Ainsi (1) entraîne que $\text{Red}_{\rho_0}(\pi_0)$ est de cardinal exactement 1. En un tel point, d'après la proposition, $n_1(s_0) = 1$ et $n_0(s_0) = 0$. Si $s_0 > 1/2$, $n_0(s_0 - 1) = n_1(s_0) = 1$ et il existe $a \in \text{Jord}_{\rho_0}(\pi_0)$ tel que $s_0 = (a+1)/2$. La remarque (2) montre qu'alors $\text{Jord}_{\rho_0}(\pi_0)$ est de la forme $\{a_{\rho_0} - 2k; k \in [0, [a_{\rho_0}/2]]\}$ et $s_0 = (a_{\rho_0}+1)/2$. L'autre possibilité est que le point de réductibilité soit $1/2$ mais alors, par (2) et l'unicité, $\text{Jord}_{\rho_0}(\pi_0) = \emptyset$.

Supposons maintenant que $\rho_0 \simeq \rho_0^*$ et que $L_{v_0}(\rho_0, r, s)$ n'a pas de pôle en $s = 0$. Alors $n_1(s_0) = 0$ pour tout $s_0 \in ]0,1[$ et $n_0(0) = n_1(1) = 1$ précisément si $\rho_0 \times \pi_0$ est irréductible. Ainsi soit $0$ est le seul point de réductibilité soit $0$ n'est pas point de réductibilité et $n_1(1) = 1$. Supposons que $n_1(1) = 1$, par définition $\text{Jord}_{\rho_0}(\pi_0) \neq \emptyset$; (1) montre qu'il n'y a qu'un point de réductibilité nécessairement dans $\mathbb{N}$ et (2) le calcule. En outre on obtient la propriété cherchée sur $\text{Jord}_{\rho_0}(\pi_0)$. Cela termine la preuve.

## 3. Paquets d'Arthur

Dans ce numéro, on ne suppose plus que la condition (C) est vérifiée, on ne suppose que (A) et (B). On impose d'abord en plus qu'il existe $\tilde{\pi}$ représentation irréductible de $\text{GL}(2n + \varepsilon, \mathbb{A})$ satisfaisant ces conditions et dont le support cuspidal satisfait à la conjecture de Ramanujan-Peterson à la place $v_0$. C'est la condition (C') de l'introduction.

Grâce aux représentations $\tilde{\pi}_{v_0} =: \tilde{\pi}_0$ où $\tilde{\pi}$ vérifie les conditions ci-dessus, on va associer des paramètres de paquets d'Arthur contenant (conjecturalement ou par défintion suivant le point de vue) $\pi_0$.

L'hypothèse mise entraîne que $\tilde{\pi}_0$ est de la forme:

$$\times_{(\sigma,b)\in K_0} \pi(\sigma, b),$$

où $\sigma$ est une série discrète d'un groupe linéaire sur $F$ de rang noté $d_\sigma$ et où $b \in \mathbb{N}$. On suppose encore que par la conjecture de Langlands locale on sache associer à chaque série discrète $\sigma$ un homomorphisme:

$$\phi_\sigma : W_F \times \mathrm{SL}(2,\mathbb{C}) \to \mathrm{GL}(d_\sigma, \mathbb{C}).$$

Pour $(\sigma, b) \in K_0$ on définit:

$$\psi_{\sigma,b} : W_F \times \mathrm{SL}(2,\mathbb{C}) \times \mathrm{SL}(2,\mathbb{C}) \to \mathrm{GL}(d_\sigma, \mathbb{C}) \times \mathrm{GL}(b, \mathbb{C}) \hookrightarrow \mathrm{GL}(bd_\sigma, \mathbb{C}),$$

de telle sorte que $\psi_{\sigma,b}$ soit le produit tensoriel de $\phi_\sigma$ et de l'homomorphisme associé à la représentation de dimension $b$ de la deuxième copie de $\mathrm{SL}(2,\mathbb{C})$.

On définit:

$$\psi_{\tilde\pi_0} := \oplus_{(\sigma,b)\in K_0} \psi_{(\sigma,b)}.$$

C'est un homomorphisme de $W_F \times \mathrm{SL}(2,\mathbb{C}) \times \mathrm{SL}(2,\mathbb{C}) \to \mathrm{GL}(2n+\varepsilon, \mathbb{C})$. Ce sont ces homomorphismes, quand $\tilde\pi$ varie soumis aux hypothèses (A), (B) et (C$'$) qui sont les bons candidats pour paramétriser les paquets d'Arthur contenant $\pi_0$; ici il n'est pas clair que leur image soit incluse dans le $L$-groupe de $G(n)$ mais cela peut se démontrer avec le théorème ci-dessous à condition de savoir interpréter l'existence de pôle à l'une des fonctions $L_{v_0}(\sigma, \wedge^2, s)$ et $L_{v_0}(\sigma, \mathrm{Sym}^2, s)$ en $s = 0$ (cf. §6, ci-dessous). On note $\Delta\psi_{\tilde\pi_0}$ l'homomorphisme de $W_F \times \mathrm{SL}(2,\mathbb{C})$ obtenu en restreignant cet homomorphisme à la diagonale de $\mathrm{SL}(2,\mathbb{C})$. Remarquons qu'à un tel homomorphisme de $W_F \times \mathrm{SL}(2,\mathbb{C})$ dans $\mathrm{GL}(2n+\varepsilon, \mathbb{C})$ on associe une décomposition en sous- $W_F \times \mathrm{SL}(2,\mathbb{C})$ représentations irréductibles de $\mathbb{C}^{2n+\varepsilon}$. Il est facile de décrire cette décomposition, en restant du côté représentations de façon à éviter l'utilisation de la conjecture de Langlands locale.

Soit $(\sigma, b) \in K_0$; comme $\sigma$ est une série discrète, d'après Zelevinski ([Z]), on sait qu'il existe $\rho_\sigma$ une représentation cuspidale d'un groupe linéaire $\mathrm{GL}(c_\sigma, F)$ et un entier $b'_\sigma$ tel que $c_\sigma b'_\sigma = d_\sigma$ et $\sigma$ soit l'unique sous-module de l'induite:

$$\times_{k \in [-(b'_\sigma-1)/2, (b'_\sigma-1)/2]} \rho_\sigma |\ |^k,$$

représentation notée usuellement $\mathrm{St}(\rho_\sigma, b'_\sigma)$; c'est ce que l'on appelle aussi une représentation de Steinberg généralisée, le cas usuel de la représentation de Steinberg étant le cas où $\rho_\sigma$ est la représentation triviale. En d'autres termes $\tilde\pi_0$ est de la forme:

(1) $$\tilde\pi_0 \simeq \times_{(\sigma,b',b) \in K_0} \pi(\mathrm{St}(\sigma, b'), b).$$

On note:

(2) $$\tilde\pi_0^A := \times_{(\sigma,b)\in K_0} \times_{\ell \in [0, \inf(b'_\sigma, b)-1]} \pi(\rho_\sigma, b'_\sigma + b - k).$$

Alors $\Delta\psi_{\tilde\pi_0}$ est naturellement associé à cette représentation: c'est un calcul élémentaire: il suffit de déterminer la représentation de $\mathrm{SL}(2,\mathbb{C})$ que l'on obtient en restreignant à la diagonale et pour cela il suffit de calculer les valeurs propres du tore diagonal intervenant dans cette représentation et on est ramené à prouver l'égalité des 2 ensembles ci-dessous, pout tout $(\sigma, b) \in K_0$:

(3) $\quad \{i+j | i \in [-(b-1)/2, (b-1)/2], j \in [-(b'_\sigma - 1)/2, (b'_\sigma - 1)/2]\}$
$$= \cup_{\ell \in [0, \inf(b'_\sigma, b)]} [-(b'_\sigma + b - \ell - 1)/2, (b'_\sigma + b - \ell - 1)/2].$$

Cela se démontre par récurrence sur $inf(b'_\sigma, b_\sigma)$.

Remarquons (voir §6 pour plus de détail) que $\tilde{\pi}_0$ est certainement de la forme (1) si $\tilde{\pi}_0$ est $G(n)$-discrète.

(C″) Dans ce qui suit on fait uniquement l'hypothèse que $\tilde{\pi}_0$ est de la forme (1). Alors $\tilde{\pi}_0^A$ est bien définie et s'interprète comme ci-dessus.

Avant d'énoncer le théorème remarquons que $\tilde{\pi}_0 \simeq \tilde{\pi}_0^A$ si $\tilde{\pi}$ vérifie la condition (C). Plus généralement,

(*) La représentation $\tilde{\pi}_0^A$ est le quotient de Langlands associé au support cuspidal de la représentation $\tilde{\pi}_0$.

On définit $\text{Jord}_{\rho_0}(\pi_0)$ comme dans l'introduction mais en utilisant $\tilde{\pi}_0^A$ au lieu de $\tilde{\pi}_0$.

THÉORÈME. *L'homomorphisme $\Delta \psi_{\tilde{\pi}_0}$ est indépendant du choix de $\tilde{\pi}$ ou plus prosaïquement $\tilde{\pi}_0^A$ est indépendant du choix de $\tilde{\pi}$ vérifiant les conditions* (A), (B), (C′) *(ou mieux* (C″)).

*Les ensembles* $\text{Red}_{\rho_0}(\pi_0)$ *pour toute représentation cuspidale unitaire irréductible d'un groupe linéaire* $\text{GL}(c, F)$ *se calculent à l'aide des ensembles* $\text{Jord}_{\rho_0}(\pi_0)$ *comme dans le théorème de l'introduction. En outre ces ensembles ont la propriété de parité et sont sans trou et sans multiplicité comme expliqué dans l'introduction.*

L'indépendance de $\tilde{\pi}_0^A$ du choix de $\tilde{\pi}$ résulte de la deuxième partie du théorème dont nous allons expliquer la démonstration.

Fixons donc $\rho_0$ et calculons $\text{Red}_{\rho_0}(\pi_0)$. On pose ici pour $s_0 \in \mathbb{R}_{\geq 0}$

$$n'_0(s_0) := -\text{ordre}_{s=s_0} L_{v_0}(\rho_0 \times \tilde{\pi}_0, s) L_{v_0}(\rho_0 \times \tilde{\pi}_0, s+1)^{-1} + \eta(s_0),$$

avec la notation $\eta(s_0)$ de 2.2 et

$$n'_1(s_0) := -\text{ordre}_{s=s_0} L_{v_0}(\rho_0 \times \tilde{\pi}_0, -s) L_{v_0}(\rho_0 \times \tilde{\pi}_0, 1-s)^{-1} + \eta(s_0) + \varepsilon_0(s_0).$$

On pose encore pour tout $s_0 \in \mathbb{R}_{\geq 0}$ $\text{irr}(s_0) = 1$ si $s_0 \in \text{Red}_{\rho_0}(\pi_0)$ et $\text{irr}(s_0) = 0$ si $s_0$ n'est pas un point de réductibilité. En 2.2, on a démontré que:

$$n'_1(s_0) = n'_0(s_0) + \text{irr}(s_0),$$

pour tout $s_0 \in \mathbb{R}_{>0}$ et $n'_0(0) + \text{irr}(0) + \varepsilon'_0 = 1$. Ceci reste vrai, puisque cela se démontre par des considérations sur les opérateurs d'entrelacement liés à l'induite $\rho_0| |^{s_0} \times \pi_0$. Ce qui n'est plus vrai est l'interprétation de $n'_1(s_0)$ et $n'_0(s_0)$ de Section 2.

On définit $n_0(s_0)$ et $n_1(s_0)$ grâce à $\tilde{\pi}_0^A$ et montrons que:

$$n'_0(s_0) - n'_1(s_0) = n_0(s_0) - n_1(s_0)$$

pour tout $s_0 \in \mathbb{R}_{>0}$. En effet le calcul des fonctions $L$ pour les paires de séries discrètes assure que:

$$\begin{aligned}\text{ordre}_{s=s_0} & L_{v_0}(\rho_0 \times \sigma, s)/L_{v_0}(\rho_0 \times \sigma, 1-s) \\ &= \text{ordre}_{s=s_0} L_{v_0}(\rho_0 \times \pi(\rho_\sigma, b'_\sigma), s)/L_{v_0}(\rho_0 \times \pi(\rho_\sigma, b'_\sigma), 1-s).\end{aligned}$$

L'égalité cherchée est alors claire sur les définitions explicites, grâce à (3) ci-dessus. En $s_0 = 0$:

$$n'_0(0) = \text{ordre}_{s=0} L_{v_0}(\rho_0 \times \tilde{\pi}_0, s)/L_{v_0}(\rho_0 \times \tilde{\pi}_0, 1-s) = n_0(0),$$

par l'argument donné ci-dessus. On termine la preuve du théorème comme celui de 2.3.

## 4. Le cas le plus général

4.1. Dans ce paragraphe, on suppose que $\pi_0$ la représentation cuspidale irréductible fixée de $G(n, F)$ est composante locale d'une situation globale satisfaisant la fonctorialité de Langlands sous la forme que lui a donnée Arthur. C'est-à-dire que l'on suppose qu'il existe une représentation cuspidale $\tilde{\pi}$ satisfaisant les hypothèses (A) et (B) de l'introduction. On explique ce que l'on peut obtenir comme résultat concernant l'irréductibilité des induites. Cela n'a d'intérêt *a priori* que technique. Toutefois on attache à $\pi_0$ une représentation irréductible $\tilde{\pi}_0^A$ de $\text{GL}(2n+\varepsilon, F)$ indépendante du choix de $\tilde{\pi}$ satisfaisant (A) et (B); c'est cette représentation et sa duale au sens de Zelevinski qui permettront la normalisation des opérateurs d'entrelacement.

On pose $\tilde{\pi}_{v_0} =: \tilde{\pi}_0$. Grâce à l'hypothèse (A), nécessairement (St(···) a été défini au début de §3):

$$\tilde{\pi}_0 \simeq \times_{(\sigma, b', b, x) \in K_0} \pi(\text{St}(\sigma, b'), b)| |^x,$$

mais où ici $\sigma$ est une représentation cuspidale unitaire, $b, b' \in \mathbb{N}$ et $x \in ]-1/2, 1/2[$. On définit: $\tilde{\pi}_0^A$ comme le quotient de Langlands de l'induite obtenue avec le support cuspidal de $\tilde{\pi}_0$ conjugué en un élément de la chambre de Weyl positive. On peut vérifier que:

$$\tilde{\pi}_0^A \simeq \times_{(\sigma, b, b', x) \in K_0} \times_{k \in [0, \inf(b, b')[} \pi(\sigma, b + b' - k)| |^x.$$

ce que l'on réécrit sous la forme:

$$\tilde{\pi}_0^A \simeq_{(\sigma, a, x) \in E_0} \pi(\sigma, a)| |^x.$$

Remarquons ici, bien que cela ne nous servira que plus loin que $E_0$ a la propriété suivante:

(*) $$(\sigma, a, x) \in E_0 \Leftrightarrow (\sigma, a, -x) \in E_0.$$

En effet, d'après les résultats de Tadic [T], on a la propriété de symétrie:

$$(\sigma, b', b, x) \in K_0 \Leftrightarrow (\sigma, b', b, -x) \in K_0$$

et (*) en résulte par construction.

On définit $\tilde{\pi}_0^L$ en prenant la conjuguée de $\tilde{\pi}_0^A$ sous l'involution de Zelevinski. Pour $\sigma$ une représentation cuspidale unitaire d'un groupe linéaire et $x \in ]-1/2, 1/2[$ on pose:

$$\operatorname{Jord}_{\sigma,x}(\pi_0) := \{a \in \mathbb{N} | (\sigma, a, x) \in E_0\}$$

ensemble compté avec multiplicité. On définit pour tout $s \in \mathbb{R}_{\geq 0}$, pour toute représentation cuspidale unitaire $\rho_0$ d'un groupe linéaire et pour tout $x \in ]-1/2, 1/2[$:

$$n_0(\rho_0, s, x) := |\{a \in \operatorname{Jord}_{\rho_0,x}(\pi_0) \mid s = x + (a-1)/2\}|.$$

Pour tout $s \in \mathbb{R}_{>1/2}$, pour tout $x \in ]-1/2, 1/2[$ et pour tout $\rho_0$ comme ci-dessus, on pose:

$$n_1(\rho_0, s, x) = |\{a \in \operatorname{Jord}_{\rho_0,x}(\pi_0) \mid s = x + (a+1)/2\}|.$$

Pour $s \in \mathbb{R}_{\geq 0}$ on pose encore:

- $\operatorname{irr}_{\rho_0}(s) = 0$ si $\rho_0 ||^s \times \pi_0$ est irréductible et
- $\operatorname{irr}_{\rho_0}(s) = 1$ si $\rho_0 ||^s \times \pi_0$ est réductible.

On définit $r_{v_0}^A(\rho_0 \times \pi_0, s)$ et $r_{v_0}^A(\rho_0^* \times \pi_0, -s)$ comme en Sections 2 et 3 et on démontre comme en loc. cit.:

4.2. PROPOSITION. *Pour tout $s_0 \in \mathbb{R}_{>1/2}$*:

$$-\operatorname{ordre}_{s=s_0}\left(r_{v_0}^A(\rho_0^* \times \pi_0, -s) r_{v_0}^A(\rho_0 \times \pi_0, s)\right)$$
$$= \sum_{x \in ]-1/2, 1/2[} n_1(\rho_0, s_0, x) - n_0(\rho_0, s_0, x) = \operatorname{irr}_{\rho_0}(s_0).$$

Le problème est de comprendre ce qui se passe pour $s_0 \in [0, 1/2]$. Aux bornes, on a la même situation qu'en Section 3 puisque la non satisfaction de la conjecture de Ramanujan-Peterson concerne les $x \neq 0$ avec $x \in ]-1/2, 1/2[$. Donc on obtient encore:

4.3. PROPOSITION. *Soit $\rho_0$ une représentation cuspidale irréductible telle que $\operatorname{Jord}_{\rho_0,0}(\pi_0) \neq \emptyset$. Alors cet ensemble est sans multiplicité sans trou* *et a la condition de parité (en particulier $\rho_0$ est autoduale). Si l'on note $a_{\rho_0,0}$ l'élément maximal de cet ensemble, alors*

$$\mathrm{Red}_{\rho_0}(\pi_0) \cap 1/2\mathbb{Z}_{\geq 0} = \{(a_{\rho_0,0}+1)/2\}.$$

*Si l'ensemble $\mathrm{Jord}_{\rho_0,0}(\pi_0) = \emptyset$ alors que $\rho_0$ est autodual, soit 0 soit 1/2 est un point de réductibilité et est l'unique point de réductibilité demi-entier positif ou nul. Le premier cas se produisant exactement s'il n'existe pas de pôle à la fonction $L_{v_0}(\rho_0, r, s)$ en $s = 0$ et le deuxième cas dans le cas opposé.*

Par contre pour les autres points, la méthode ne donne aucun résultat mais montrons toutefois:

4.4. PROPOSITION. *Les représentations $\tilde{\pi}_0^A$ et $\tilde{\pi}_0^L$ sont uniquement déterminées par $\pi_0$, i.e. sont indépendantes du choix de $\tilde{\pi}$.*

Ce n'est pas complètement évident: le problème est de déterminer pour tout $x \in ]-1/2, 1/2[$ et pour tout $\sigma$, $\mathrm{Jord}_{\sigma,x}(\pi_0)$. On fixe $\sigma$ que je préfère appeler $\rho_0$ pour être conforme aux notations précédentes. Le cas de $x = 0$ est traité par la proposition ci-dessus. Supposons donc que $x \neq 0$; on utilise le fait que:

(1) $$\mathrm{Jord}_{\rho_0,x}(\pi_0) = \mathrm{Jord}_{\rho_0,-x}(\pi_0).$$

Cela résulte des propriétés des représentations unitaires de $\mathrm{GL}(n)$ démontrées par Tadic [T] (cf. (*) ci-dessus). On procède alors ainsi: l'ensemble des points de réductibilité, non 1/2 entiers, coïncide avec l'ensemble suivant obtenu comme différence d'ensembles considérés avec multiplicité:

(2) $\quad E := \{x + (a+1)/2; x \in ]-1/2, 1/2[-\{0\}, a \in \mathrm{Jord}_{\rho_0,x}(\pi_0)\}$,

(2') $\quad E := -\{x + (a-1)/2; x \in ]-1/2, 1/2[-\{0\}, a \in \mathrm{Jord}_{\rho_0,x}(\pi_0)\}$.

En outre grâce à 4.2, on sait que l'ensemble $E$ est sans multiplicité. Montrons que $E$ détermine l'ensemble $E'$

$$\{(x,a) | x \in ]0, 1/2[; a \in \mathrm{Jord}_{\rho_0,x}(\pi_0)\}.$$

En tenant compte de (1) cela suffira. On met sur $E'$ l'ordre suivant

$$(x,a) \geq (x',a') \Leftrightarrow x + (a+1)/2 \geq x' + (a'+1)/2.$$

Remarquons que la positivité mise dans la définition de $E'$ entraîne que si $x + (a+1)/2 = x' + (a'+1)/2$ alors $(x,a) = (x',a')$ au sens usuel. On détermine les éléments de $E'$ avec leur multiplicité par récurrence décroissante. On fixe $y \in \mathbb{R}_{>0}$ et on suppose que l'on a déterminé (avec leur multiplicité) les éléments $(x,a)$ de $E'$ tel que $x + (a+1)/2 > y$ et on va déterminer les éléments de $E'$ maximaux avec la propriété $x + (a+1)/2 \leq y$. On peut évidemment supposer que la valeur $y$ est atteinte. et on va donc déterminer les éléments de $E'$ tels que $x + (a+1)/2 = y$. D'abord on remarque que $y$ détermine uniquement $x$ que l'on note donc $x_y$. De plus $y$ est dans l'ensemble écrit en (2). On note $m_+(y)$ la multiplicité avec laquelle $y$ apparaît dans cet ensemble (2) définissant $E$ et $m_-(y)$ celle avec laquelle il apparaît dans l'ensemble (2'). Grâce à la connaissance de $E$, on connait:

$$m_+(y) - m_-(y).$$

On montre d'abord que cela suffit pour déterminer $m_+(y)$: en effet $m_-(y)$ est déterminé par l'ensemble des éléments $(x,a)$ tels que $x \in\ ]-1/2, 1/2[$, $a \in$ Jord$_{\rho_0,x}(\pi_0)$ tels que $x + (a-1)/2 = y$. Or un tel couple vérifie certainement

$$x + (a+1)/2 > y$$

et encore *a fortiori* $|x| + (a+1)/2 > y$. En tenant compte de (1) et de l'hypothèse de récurrence ces couples sont parfaitement connus ainsi que leur multiplicité, d'où la connaissance de $m_+(y)$. Les éléments qui contribuent à la multiplicité $m_+(y)$ sont d'une part le couple $(x_y, 2y - x_y - 1)$ et d'autre part les éléments $(x,a)$ vérifiant

$$x + (a+1)/2 = y, \quad x \in\ ]-1/2, 0[.$$

Ces éléments du deuxième type sont connus et cela permet de déterminer la multiplicité de $(x_y, 2y - x_y - 1)$.

Une fois que $m_+(y)$ est connu, on connait tous les ensembles Jord$_{\rho_0,x}(\pi_0)$ pour $x \in\ ]0, 1/2[$ et il ne reste plus qu'à utiliser encore (1) pour conclure.

## 5. Normalisation des opérateurs d'entrelacement

Dans ce paragraphe, on suppose que $\pi_0$ la représentation cuspidale irréductible fixée de $G(n,F)$ satisfait les hypothèses relativement générales de Section 4. En Section 4 spécialement 4.4, on a associé à $\pi_0$ une représentation $\tilde{\pi}_0^A$ de GL$(2n+\varepsilon, F)$. On note $\tilde{\pi}_0^L$ la représentation de GL$(2n+\varepsilon, F)$ duale pour l'involution de Zelevinski de $\tilde{\pi}_0^A$; cf. 4.1.

On pose:

$r_{v_0}^A(\rho_0 \times \pi_0, s)$
$\quad := L_{v_0}(\rho_0 \times \tilde{\pi}_0^A, s) L_{v_0}(\rho_0, r, 2s) L_{v_0}(\rho_0 \times \tilde{\pi}_0^A, s+1)^{-1} L_{v_0}(\rho_0, r, 2s+1)^{-1}$

et

$r_{v_0}^L(\rho_0 \times \pi_0, s)$
$\quad := L_{v_0}(\rho_0 \times \tilde{\pi}_0^L, s) L_{v_0}(\rho_0, r, 2s) L_{v_0}(\rho_0 \times \tilde{\pi}_0^L, s+1)^{-1} L_{v_0}(\rho_0, r, 2s+1)^{-1}.$

On normalise alors les opérateurs d'entrelacements des deux façons suivantes:

$$N_{v_0}^A(\rho_0 \times \pi_0, s) := r_{v_0}^A(\rho_0 \times \pi_0, s)^{-1} M_{v_0}(\rho_0 \times \pi_0, s)$$

et
$$N_{v_0}^L(\rho_0 \times \pi_0, s) := r_{v_0}^L(\rho_0 \times \pi_0, s)^{-1} M_{v_0}(\rho_0 \times \pi_0, s).$$

On peut évidemment donner des définitions analogues en remplaçant $\rho_0$ par $\rho_0^*$.

THÉORÈME.　(i) *A une fonction holomorphe inversible près*:
$$r_{v_0}^A(\rho_0^* \times \pi_0, -s) r_{v_0}^A(\rho_0 \times \pi_0, s) = r_{v_0}^L(\rho_0^* \times \pi_0, -s) r_{v_0}^L(\rho_0 \times \pi_0, s);$$
*et cette fonction diffère de*
$$M_{v_0}(\rho_0^* \times \pi_0, -s) M_{v_0}(\rho_0 \times \pi_0, s)$$
*par une fonction holomorphe inversible.*

(ii) *Les opérateurs $N_{v_0}^A(\rho_0 \times \pi_0, s)$ et $N_{v_0}^L(\rho_0 \times \pi_0, s)$ sont holomorphiquement définis pour tout $s$ avec $\operatorname{Re} s \geq 0$; le deuxième n'a pas de zéro dans cet ensemble et le premier en a précisément aux points $x_0 := s_0 + iy_0$ où $s_0$ vérifie les conditions suivantes*:

$s_0 \neq 0$,

*il existe $x \in\ ]-1/2, 1/2[\ $ et $a \in \operatorname{Jord}_{\rho_0| \ |^{iy_0}, x}(\pi_0)$ tels que $s_0 = x + (a-1)/2$.*

(iii) *Les produits*
$$N_{v_0}^A(\rho_0^* \times \pi_0, -s) N_{v_0}^A(\rho_0 \times \pi_0, s)$$
*et*
$$N_{v_0}^L(\rho_0^* \times \pi_0, -s) N_{v_0}^L(\rho_0 \times \pi_0, s)$$
*sont des fonctions holomorphes inversibles.*

L'égalité des ordres des 2 fonctions en (i) a déjà été vue en Section 3. Quitte à tensoriser $\rho_0$ par un caractère du type $|\ |^{iy_0}$ où $y_0 \in \mathbb{R}$, on se ramène à faire l'étude au voisinage des points de $\mathbb{R}$. D'après Harish-Chandra, pour tout $s \in \mathbb{R} - \{0\}$ le produit des opérateurs d'entrelacements non normalisés est holomorphe avec des zéros simples en $s_0$ tel que l'induite soit réductible. En 4.2 on a précisément démontré la même assertion pour le produit:
$$r_{v_0}^A(\rho_0^* \times \pi_0, -s) r_{v_0}^A(\rho_0 \times \pi_0, s)$$
du moment que $s > 0$. Supposons maintenant que $s < 0$; alors on applique ce que l'on vient de voir avec l'induite $\rho_0^*|\ |^{-s} \times \pi_0$ dont l'irréductibilité est équivalent à l'irréductibilité de $\rho_0|\ |^s \times \pi_0$. Et on obtient que le produit:
$$r_{v_0}^A(\rho_0 \times \pi_0, s) r_{v_0}^A(\rho_0^* \times \pi_0, -s)$$
est nul à l'ordre 1 exactement au point de réductibilité de cette induite si $-s > 0$, d'où (i).

(i) entraîne (iii).

Démontrons (ii) sans supposer que $s \in \mathbb{R}$ mais seulement $\operatorname{Re} s \geq 0$. On calcule explicitement avec les notations reprises de 4.1:

$$\tilde{\pi}_0^A =: \times_{(\sigma,a,x) \in E_0} \pi(\sigma, a),$$

$$r_{v_0}^A(\rho_0 \times \pi_0, s) = \times_{(\sigma,a,x) \in E_0} L_{v_0}(\rho_0 \times \sigma^*, s - (a-1)/2 - x) L_{v_0}(\rho_0, r, 2s)$$
$$\cdot L_{v_0}(\rho_0 \times \sigma^*, s + (a-1)/2 - x + 1)^{-1} L(\rho_0, r, 2s + 1)^{-1}.$$

(J'ai utilisé l'autodualité de $\tilde{\pi}_0^A$ qui résulte de celle de $\tilde{\pi}_{v_0}$; c'est uniquement pour avoir des formules plus jolies.) Les dénominateurs sont holomorphes pour tout $s$ tel que $\operatorname{Re} s \geq 0$. Les numérateurs introduisent éventuellement des pôles. Soit $x_0$ un tel pôle; on écrit $x_0 = s_0 + iy_0$ avec $s_0, y_0 \in \mathbb{R}$. Les facteurs qui donnent des pôles doivent satisfaire:

$$\rho_0| \ |^{iy_0} \simeq \sigma$$

et $s_0 = (a-1)/2 + x$. Comme l'opérateur d'entrelacement non normalisé n'a pas de pôles en dehors de l'axe unitaire, ces pôles introduisent des 0 pour les opérateurs d'entrelacement normalisés grâce à $\tilde{\pi}_0^A$.

En $s_0 = 0$: on a des pôles soit si $L_{v_0}(\rho_0| \ |^{iy_0}, r, s)$ a un pôle en $s = 0$ soit s'il existe $(\sigma, 0, 1) \in E_0$ tel que $\rho_0| \ |^{iy_0} \simeq \sigma$. En remplaçant $\rho_0$ par $\rho_0| \ |^{iy_0}$ on a vu que ces deux possibilités s'excluaient l'une l'autre et l'une des ces deux propriétés est réalisée précisément quand 0 n'est pas point de réductibilité pour l'induite $\rho_0| \ |^{iy_0} \times \pi_0$. En d'autres termes $r_{v_0}^A(\rho_0 \times \pi_0, s)$ a un pôle en $iy_0$ si et seulement si $M_{v_0}(\rho_0 \times \pi_0, s)$ a un pôle en $iy_0$; dans ce cas les pôles sont simples. Cela prouve l'holomorphie de $N_{v_0}^A(\rho_0 \times \pi_0, s)$ en tout $s$ tel que $\operatorname{Re} s \geq 0$. En outre le facteur de normalisation introduit des zéros précisément aux points donnés dans l'énoncé.

Il faut faire un raisonnement du même type en utilisant $\tilde{\pi}_0^L$. On a:

$$\tilde{\pi}_0^L = \times_{(\sigma,x,a) \in E_0} \operatorname{St}(\sigma, a)| \ |^x,$$

où St est la représentation de Steinberg généralisée. On calcule:

$$r_{v_0}^L(\rho_0 \times \pi_0, s) = \times_{(\sigma,a) \in E_0} L_{v_0}(\rho_0 \times \sigma, s + (a-1)/2 - x) L_{v_0}(\rho_0, r, 2s)$$
$$\cdot L_{v_0}(\rho_0 \times \sigma, s + (a-1)/2 + 1 - x)^{-1} L(\rho_0, r, 2s + 1)^{-1}.$$

Les dénominateurs n'introduisent donc pas de zéros. Les numérateurs sont holomorphes en tout $s$ tel que $\operatorname{Re} s \geq 0$ sauf éventuellement en $s = iy_0$ imaginaire. Et cela se produit exactement si $\operatorname{Jord}_{(\rho_0| \ |^{iy_0}, 0)}(\pi_0)$ contient 1. On termine la preuve exactement comme ci-dessus.

## 6. Correspondance de Langlands

Reprenons les hypothèses de Section 4, c'est-à-dire que l'on fixe $\pi_0$ une représentation cuspidale irréductible de $G(n, F)$ et que l'on suppose que les hypothèses (A) et (B) de l'introduction sont satisfaites. En Section 4, on a associé à $\pi_0$ la représentation $\tilde{\pi}_0^L$. Elle est de la forme:

$$\tilde{\pi}_0^L \simeq \times_{(\sigma,a,x)\in K_0} \mathrm{St}(\sigma, a)|\ |^x$$

avec la propriété que si $(\sigma, a, x) \in K_0$ alors $(\sigma, a, -x)$ est aussi dans $K_0$ avec la même multiplicité.

[B-H-K] ont asssocié à toute représentation cuspidale irréductible, $\sigma$, d'un groupe linéaire $\mathrm{GL}(d_\sigma, F)$ un morphisme:

$$\psi_\sigma : W_F \to \mathrm{GL}(d_\sigma, \mathbb{C})$$

qui est l'unique candidat possible pour satisfaire aux conjectures de Langlands; ce genre de correspondance se trouve déjà dans des papiers antérieurs de Henniart. La conjecture de Langlands est maintenant démontrée; cf. l'introduction. Je vais supposer en plus ici que l'application $\sigma \mapsto \psi_\sigma$ a la propriété suivante:

- Si $L_{v_0}(\sigma, \wedge^2 \mathbb{C}^{d_\sigma}, s)$ a un pôle en $s = 0$, alors $d_\sigma$ est paire et $\psi_\sigma$ se factorise par le groupe symplectique $\mathrm{Sp}(d_\sigma, \mathbb{C})$ (à conjugaison près évidemment);
- Si $L_{v_0}(\sigma, \mathrm{Sym}^2 \mathbb{C}^{d_\sigma}, s)$ a un pôle en $s = 0$, alors $\psi_\sigma$ se factorise (à conjugaison près) par le groupe orthogonal $\mathrm{O}(d_\sigma, \mathbb{C})$.

A une série discrète du type $\mathrm{St}(\sigma, a)$, on associe alors le morphisme:

$$\psi_{\sigma,a} : W_F \times \mathrm{SL}(2, \mathbb{C}) \to \mathrm{GL}(d_\sigma, \mathbb{C}) \times \mathrm{GL}(a, \mathbb{C}) \hookrightarrow \mathrm{GL}(\mathbb{C}^{d_\sigma} \otimes \mathbb{C}^a) = \mathrm{GL}(ad_\sigma, \mathbb{C}).$$

Pour $x \in \mathbb{R}$, on définit aussi le morphisme produit du morphisme $\psi_\sigma$ et du morphisme correspondant à la représentation de dimension $a$ de $\mathrm{SL}(2, \mathbb{C})$:

$$\lambda_x : W_F \to \mathbb{C}^*$$

$$x \in W_F \hookrightarrow |x| \in \mathbb{C}^*.$$

On note $\psi_{\sigma,a,x}$ le morphisme obtenu par multiplication de $\psi_{\sigma,a}$ avec $\lambda_x$ vu comme morphisme de $W_F \times \mathrm{SL}(2, \mathbb{C})$ dans le centre de $\mathrm{GL}(ad_\sigma, \mathbb{C})$ trivial sur $\mathrm{SL}(2, \mathbb{C})$.

On pose:

$$\psi_{\pi_0} : W_F \times \mathrm{SL}(2, \mathbb{C}) \to \mathrm{GL}(2n + \varepsilon, \mathbb{C})$$

$$\psi_{\pi_0} = \oplus_{(\sigma,a,x)\in K_0} \psi_{\sigma,a,x}.$$

Chacun des morphismes $\psi_{\sigma,a,x}$ n'a de sens qu'à conjugaison près et donc $\psi_{\pi_0}$ n'est défini qu'à conjugaison près.

On fixe une inclusion de ${}^L G(n) \hookrightarrow \mathrm{GL}(2n + \varepsilon, \mathbb{C}) \times W_F$: pour cela, on rappelle (cf. [M]) que:

- Si $G$ est un groupe symplectique, ${}^L G(n) \simeq \mathrm{SO}(2n+1, \mathbb{C}) \times W_F$;
- Si $G$ est un groupe spécial orthogonal d'une forme de dimension impaire (déployée ou non) ${}^L G(n) \simeq \mathrm{Sp}(2n) \times W_F$; on a expliqué en [M] comment, dans le yoga des paquets d'Arthur, se distingue le cas quasi-déployé du cas non quasi-déployé; cela se fait non pas sur le paramètre

$$W_F \times \mathrm{SL}(2,\mathbb{C}) \times \mathrm{SL}(2,\mathbb{C}) \to {}^L G(n) \times W_F$$

   mais sur les caractères du centralisateur de ce paramètre qui interviennent dans l'étude plus fine des paquets d'Arthur. Cela n'interviendra donc pas ici, puisque nos arguments sont insuffisants pour déterminer ces caractères. En [M], on avait considéré le cas du groupe orthogonal et suivant [Ad], on avait ajouté un facteur $\mathbb{Z}/2\mathbb{Z}$ qui n'est pas à mettre ici puisque l'on a travaillé avec la composante neutre.

- Si $G$ est un groupe orthogonal d'une forme de dimension paire, nous prenons comme définition de ${}^L G(n) :\simeq \mathrm{O}(2n, \mathbb{C}) \times W_F$ (cela est inspiré aussi de [Ad]). La distinction du cas déployé et du cas non déployé comme ci-dessus n'intervient pas ici.

Maintenant que l'on a écrit ${}^L G(n)$ comme produit direct d'un groupe réductif complexe avec $W_F$, on oublie ce facteur $W_F$ en réduisant l'inclusion ci-dessus en une inclusion:

$$ {}^L G(n) \hookrightarrow \mathrm{GL}(2n+\varepsilon, \mathbb{C}). $$

PROPOSITION. *Avec les hypothèses et notations ci-dessus, à conjugaison près $\psi_{\pi_0}$ a son image incluse dans ${}^L G(n)$.*

Les arguments développés ci-dessous sont similaires à ceux de [$S_3$] et [$S_4$].

Fixons d'abord $x \in ]0, 1/2[$ et $\sigma$ une représentation cuspidale unitaire tels que $\mathrm{Jord}_{\sigma,x} \neq \emptyset$; pour tout $a$ dans cet ensemble, on conjugue $\psi_{\sigma,x,a} \oplus \psi_{\sigma,-x,a}$ de sorte que son image soit incluse dans $\mathrm{GL}(2n+\varepsilon, \mathbb{C})$ avec l'image de $\psi_{\sigma,x,a}$ incluse dans le groupe des automorphismes d'un espace isotrope de dimension $ad_\sigma$. Regardons maintenant le cas où $x=0$: fixons $\sigma$ telle que $\mathrm{Jord}_{\sigma,0} \neq \emptyset$ et fixons $a$ dans cet ensemble. En 4.3, on a montré que $a$ est pair si $L_{v_0}(\sigma,r,s)$ a un pôle en $s=0$ et est impair sinon. On remarque que $r \simeq \wedge^2 \mathbb{C}^{d_\sigma}$ si ${}^L G(n)$ est un groupe orthogonal et $r \simeq \mathrm{Sym}^2 \mathbb{C}^{d_\sigma}$ si ${}^L G(n)$ est un groupe symplectique. Ainsi les hypothèses faites sur l'interprétation des pôles de $L_{v_0}(\sigma,r,s)$, disent que s'il y a un pôle à $L_{v_0}(\sigma,r,s)$ en $s=0$, alors $\psi_\sigma$ est d'image un groupe d'automorphismes symplectiques (resp. orthogonaux) si ${}^L G(n)$ est un groupe orthogonal (resp. symplectique). Or une représentation de dimension paire de $\mathrm{SL}(2,\mathbb{C})$ est symplectique. Nous avons ainsi montré que si $L_{v_0}(\sigma,r,s)$ a un pôle en $s=0$ alors l'image de $\psi_{\sigma,0,a}$ est formé d'automorphismes symplectique si ${}^L G(n)$ est symplectique et orthogonaux sinon.

Supposons maintenant que $L_{v_0}(\sigma, r, s)$ n'a pas de pôle en $s = 0$. On a montré que $\sigma \simeq \sigma^*$ si $\mathrm{Jord}_{\sigma,0} \neq \emptyset$ et donc $L_{v_0}(\sigma \times \sigma, s)$ a un pôle en $s = 0$. Ainsi si $^L G(n)$ est symplectique $L_{v_0}(\sigma, \wedge^2 \mathbb{C}^{d_\sigma}, s)$ a un pôle en $s = 0$; alors $\psi_\sigma$ est formé d'automorphismes symplectiques et tout élément, $a$, de $\mathrm{Jord}_{\sigma,0}$ est impair. Une représentation de dimension impaire de $\mathrm{SL}(2,\mathbb{C})$ est orthogonale et on a encore le résultat que si $^L G(n)$ est symplectique et si $L_{v_0}(\sigma, r, s)$ n'a pas de pôle en $s = 0$, $\psi_{\sigma,0,a}$ est formé d'automorphismes symplectiques. On a un raisonnement identique dans le cas où $^L G(n)$ est orthogonal. Cela termine la preuve.

*Remarque.* Avec les notations et hypothèses ci-dessus, le morphisme

$$\psi_{\pi_0} : W_F \times \mathrm{SL}(2,\mathbb{C}) \to{}^L G(n)$$

a son image incluse dans un parabolique propre de $^L G(n)$ si et seulement si il existe $x \in \,]0, 1/2[$ et $\sigma$ une représentation cuspidale unitaire telle que $\mathrm{Jord}_{\sigma,x} \neq \emptyset$.

Institut de Mathématiques de Jussieu, Université Paris 7, Paris, France
*E-mail address*: moeglin@math.jussieu.fr

[S$_4$]   F. SHAHIDI, Twisted endoscopy and reducibility of induced representations for $p$-adic groups, *Duke Math. J.* **66** (1992), 1–41.

[Si]    A. SILBERGER, Special representations of reductive $p$-adic groups are not integrable, *Ann. of Math.* **111** (1980), 571–587.

[T]     M. TADIC, Classification of unitary representations in irreducible representations of the general linear group (non-Archimedean case), *Ann. Sci. École Norm. Sup.* **19** (1986), 335–382.

[Z]     A. ZELEVINSKY, Induced representations of reductive $p$-adic groups II. On irreducible representations of GL($n$), *Ann. Sci. École Norm. Sup.* **13** (1980), 165–210.

[Zh]    Y. ZHANG, $L$-packets and reducibilities, *J. Reine Angew. Math.* **510** (1999), 83–102.